\documentclass{article}

\usepackage{amsfonts}
\usepackage{mathrsfs}
\usepackage{dsfont}
\usepackage{amssymb}
\usepackage{amsthm}
\usepackage{theoremref}

\newtheorem{thm}{Theorem}[section]
\newtheorem{lem}[thm]{Lemma}
\newtheorem{cor}[thm]{Corollary}
\newtheorem{prp}[thm]{Proposition}
\newtheorem{qst}[thm]{Question}

\theoremstyle{definition}
\newtheorem{dfn}[thm]{Definition}

\usepackage{a4wide}

\title{
\textbf{\huge{Filters in C$^*$-Algebras}}\\
}
\author{{\Large Tristan Bice}}

\begin{document}

\maketitle

\paragraph{Abstract} In this paper we analyze states on C$^*$-algebras and their relationship to filter-like structures of projections and positive elements in the unit ball.  After developing the basic theory we use this to investigate the Kadison-Singer conjecture, proving its equivalence to an apparently quite weak paving conjecture and the existence of unique maximal centred extensions of projections coming from ultrafilters on $\omega$.  We then prove that Reid's positive answer to this for q-points in fact also holds for rapid p-points, and that maximal centred filters are obtained in this case.  We then show that consistently such maximal centred filters do not exist at all meaning that, for every pure state $\phi$ on the Calkin algebra, there exist projections $p$ and $q$ such that $\phi(p)=1=\phi(q)$, even though $\phi(r)$, for projections $r\leq p,q$, is bounded strictly below $1$.  Lastly we investigate towers, using cardinal invariant equalities to construct towers on $\omega$ that do and do not remain towers when canonically embedded into the Calkin algebra.  Finally we show that consistently all towers on $\omega$ remain towers under this embedding.\footnote{This research has been supported by a Japanese Government Mombukagakusho (Ministry of Education, Culture, Sports, Science and Technology) Scholarship.}\\

\section{Introduction}

States on C$^*$-algebras and their relation to other objects in C$^*$-algebras, like closed left ideals and closed hereditary cones, have been studied for quite some time.  Despite this, some basic questions about states remain unresolved, like the long-standing Kadison-Singer conjecture.  This paper aims to provide another, more order theoretic, perspective on states by investigating their relation to certain filter-like objects.  We also demonstrate how analyzing these objects can be used to investigate these outstanding problems and also give some new unexpected results.  It is our hope that this is just the beginning, and further research in this direction will prove even more fruitful.

In \S\ref{states} we present the basic theory of states and their relation to what we have called norm centred sets.  These first appeared under the name `quantum filters' in some joint unpublished work of Farah and Weaver, and were further developed in Farah's unpublished notes \cite{e}.  With the exception of \thref{Pperp}, the results of this section are originally from \cite{e}, although the proofs and presentation differ somewhat.  The author would like to thank Ilijas Farah for making these notes available and allowing the relevant results to be reproduced here, as well as for providing insightful comments on earlier versions of this paper.

We next define norm filters and prove some of their basic properties in \S\ref{nfs}.  The main purpose of this section is merely to demonstrate that these norm filters appear to be the natural analogs of filters in the general C$^*$-algebra context.

In \S\ref{op} we restrict our attention to C$^*$-algebras $A$ of real rank zero, where it suffices to look at norm centred subsets of orthogonal projections in $A$, rather than arbitrary positive elements in the unit ball.  Further restricting our attention to the case when the canonical order on these projections is countably (downwards) closed, we show how norm centred sets correspond to centred sets in the usual order theoretic sense.  We then prove some basic order theoretic properties of centred subsets of projections in this case.

Next, in \thref{KS}, we show how this fairly elementary theory can be used to obtain a paving conjecture equivalent of the Kadison-Singer conjecture which at first appears to be siginificantly weaker than the conjecture.  In the same theorem, we also show how the Kadison-Singer conjecture is equivalent to the simple statement that every ultrafilter on the natural numbers has a unique maximal centred extension when mapped canonically to the Calkin algebra.  With different terminology, Reid showed that this holds for q-points (see \thref{Reid}) and in \cite{u} Problem 2 it was asked if this also holds for p-points.  In \thref{rapp} we show that this does at least hold for rapid p-points.  Moreover, in this case the maximal centred extension is a filter, and hence an ultrafilter.  We then investigate ultrafitler extensions in the Calkin algebra, showing in \thref{qpufb} that q-points also give rise to unique ultrafilter extensions, although they may differ from their unique maximal centred extensions, by \thref{nonp}.

In \S\ref{mcf} we investigate maximal centred filters of projections in the Calkin algebra, i.e. ultrafilters that are maximal not just among all proper filters, but among all centred sets.  Our main result, \thref{nupf}, is that consistently they do not exist at all, specifically that this holds in a well known model of ZFC without p-points.  When translated back into the language of states, this yields the somewhat surprising result given in \thref{nupfcor}, namely that it is consistent with ZFC that, for every pure state $\phi$ on the Calkin algebra, there exist projections $p$ and $q$ such that $\phi(p)=1=\phi(q)$ even though the set of values $\phi(r)$, for projections $r\leq p,q$, has an upper bound strictly below $1$.

Lastly, in \S\ref{towers}, we investigate towers of projections in the Calkin algebra, specifically those arising from towers of subsets of $\omega$.  Despite the fact that towers are a special case of filters, this section does not use any results from previous sections (except for \thref{countdirect} which is also proved independently of other results in this paper) and may consequently be read in isolation.  It does, however, require knowledge of some cardinal invariants of the continuum, as well as forcing (which is also required for \S\ref{mcf}).  This is because we use cardinal invariant equalities to construct towers that do and do not remain towers when canonically embedded in the Calkin algebra, in \thref{b=t} and \thref{nonM=m} respectively.  Finally, we use an unpublished result of Brendle's from \cite{q} to prove in \thref{MAsigma} that iterating with all $\sigma$-centred forcings yields a model where all towers on $\omega$ remain towers under the canonical embedding into the Calkin algebra.

\section{States and Norm Centred Sets}\label{states}

First, let us set out some notation.  For a subset $A$ of a C$^*$-algebra, $A^1$ denotes the elements of the unit ball in $A$, $A_+$ denotes the positive elements in $A$ and $\mathcal{P}(A)$ denotes the orthogonal projections in $A$.  We let $\mathbb{S}(A)$ denote the states of $A$, i.e. \[\mathbb{S}(A)=\{\phi\in A^*:||\phi||=1\textrm{ and }\phi[A_+]\subseteq\mathbb{R}_+\},\] while $\mathbb{P}(A)$ denotes the pure states of $A$, i.e. the extreme points of $\mathbb{S}(A)$.  Given $\mathcal{A}\subseteq A^1_+$, we let $\mathbb{S}(\mathcal{A})=\{\phi\in\mathbb{S}(A):\forall a\in\mathcal{A}(\phi(a)=1)\}$ and $\mathbb{P}(\mathcal{A})=\mathbb{S}(\mathcal{A})\cap\mathbb{P}(A)$.  Also, given $\phi\in\mathbb{S}(A)$, we let $\phi^1_+=\{a\in A^1_+:\phi(a)=1\}$ and $\mathcal{P}(\phi)=\{a\in\mathcal{P}(A):\phi(a)=1\}$.

Throughout this section we will make use of the well known Gelfand-Naimark-Segal construction, namely that, for all $\phi\in\mathbb{S}(A)$, there exists a representation $\pi_\phi$ of $A$ on a Hilbert space $H_\phi$ and (cyclic) $v_\phi\in H^1$ with $\phi(a)=\langle\pi_\phi(a)v_\phi,v_\phi\rangle$, for $a\in A$.  We do not go into its proof, suffice to say that it uses the Cauchy-Schwarz inequality.  The proofs of the corresponding theorems in \cite{e} use the Cauchy-Schwarz inequality directly, rather than the GNS construction as done here.

We will also use the following elementary results.  Firstly, whenever $\phi\in\mathbb{S}(A)$, $a\in A$ and $b\in\phi^1_+$, we have $\phi(ab)=\phi(a)$.  To see this, simply note that $b\in\phi^1_+$ means that $\langle\pi_\phi(b) v_\phi,v_\phi\rangle=1$ which in turn yields $\pi_\phi(b)v_\phi=v_\phi$ and hence $\phi(ab)=\langle\pi_\phi(a)\pi_\phi(b)v_\phi,v_\phi\rangle=\langle\pi_\phi(a)v_\phi,v_\phi\rangle=\phi(a)$.  On the other hand, whenever $\phi\in\mathbb{S}(A)$, $a\in A^1$ and $b\in A^1_+$ and $|\phi(ab)|=1$ then $\phi(b)=1$.  For if we had $\phi(b)<1$ then that would mean $||\pi_\phi(b)v_\phi||<1$ and hence $|\phi(ab)|\leq||\pi_\phi(a)\pi_\phi(b)v_\phi||\leq||\pi_\phi(b)v_\phi||<1$, a contradiction.  This, plus the fact that states have norm one, implies that $\phi^1_+$ is always norm centred according to the following definition.

\begin{dfn}
Take a C$^*$-algebra $A$.  We call $\mathcal{A}\subseteq A^1_+$ \emph{norm centred} if $||a_1\ldots a_n||=1$ for all $n\in\omega$ and $a_1,\ldots,a_n\in\mathcal{A}$.
\end{dfn}

Our terminology differs slightly from that in \cite{e}, where norm centred sets are called quantum filters (and the definition in \cite{e} refers only to subsets of non-zero projections, rather than arbitrary elements of $A^1_+$ as done here).  We believe that it is rather the norm filters (see \thref{normfilter}) that constitute the natural quantum analog of a filter, as we discuss later on in \S\ref{nfs}.

As just mentioned, states give rise to norm centred sets, and we now show that, conversely, norm centred sets give rise to states.

\begin{thm}\thlabel{3.1(4)}
Assume $A$ is a C$^*$-algebra, $\mathcal{A}\subseteq A^1_+$ is norm centred, $b\in A^1_+$ and \[\lambda=\inf\{||a_1\ldots a_nba_n\ldots a_1||:n\in\omega\textrm{ and }a_1,\ldots,a_n\in\mathcal{A}\}.\]  Then there exists $\phi\in\mathbb{P}(A)$ such that $\phi(a)=1$, for all $a\in\mathcal{A}$, and $\phi(b)=\lambda$.
\end{thm}

\paragraph{Proof:} Let us first assume that $\mathcal{A}=\{a\}$, for some $a\in A^1_+$.  For each $\epsilon\in(0,1)$ let \[\mathbb{S}_\epsilon=\{\phi\in\mathbb{S}(A):\phi(a)\geq1-\epsilon\textrm{ and }\phi(b)\geq\lambda-\epsilon\}.\]  We claim that each $\mathbb{S}_\epsilon$ is not empty.  To see this, let $\delta>0$ be such that $1-\lambda/(\lambda+\delta)\leq\epsilon^2/8$ and let $n\in\omega$ be such that $||a^nba^n||\leq\lambda+\delta$, and hence $a^{2n}ba^{2n}\leq(\lambda+\delta)a^{2n}$.  By \cite{a} Theorem 1.7.2, we have $\phi\in\mathbb{S}(A)$ be such that $\phi(a^{2n}ba^{2n})=||a^{2n}ba^{2n}||$.  But then $\lambda\leq\phi(a^{2n}ba^{2n})\leq(\lambda+\delta)\phi(a^{2n})$, i.e. $\phi(a^{2n})\geq\lambda/(\lambda+\delta)\geq1-\epsilon^2/8\geq1-\epsilon$.  Thus \[||v-\pi(a^{2n})v||=\sqrt{1-\langle v,\pi(a^{2n})v\rangle-\langle\pi(a^{2n})v,v\rangle+||\pi(a^{2n})v||^2}\leq\sqrt{2-2\lambda/(\lambda+\delta)}\leq\epsilon/2.\]  As $\phi(a^{2n}ba^{2n})=\langle\pi(b)\pi(a^{2n})v,\pi(a^{2n})v\rangle$ and $\phi(b)=\langle\pi(b)v,v\rangle$, it follows that \[|\phi(b)-\phi(a^{2n}ba^{2n})|\leq2||v-\pi(a^{2n})v||\leq\epsilon,\] and hence $\phi(b)\geq\lambda-\epsilon$.

The claim is thus proved and hence $\mathbb{S}_\epsilon$, for $\epsilon>0$, is a collection of non-empty subsets of $A^*$ with the finite intersection property (as $\delta<\epsilon\Rightarrow\mathbb{S}_\delta\subseteq\mathbb{S}_\epsilon$).  For all $\epsilon>0$, $\mathbb{S}_\epsilon$ is closed in the weak$^*$-topology of $A^*$ (so long as $A$ is unital, otherwise look at the states on its unitization) and contained in $A^{*1}$, and hence compact by the Banach-Alagolu Theorem.  Thus $\bigcap_{\epsilon>0}\mathbb{S}_\epsilon$ is non-empty and any $\phi\in\bigcap_{\epsilon>0}\mathbb{S}_\epsilon$ will satisfy $\phi(a)=1$ and $\phi(b)=\lambda$.

For the general case of non-singleton $\mathcal{A}$, take any $a_1,\ldots,a_n\in\mathcal{A}$ and let $a=a_n\ldots a_2a_1a_2\ldots a_n$.  By the singleton case just proved, there exists $\phi\in\mathbb{S}(A)$ such that $\phi(a)=1$ and $\phi(b)=\inf||a^nba^n||\geq\lambda$.  But $\phi(a)=1$ implies that $\phi(a_k)=1$, for $k=1,\ldots,n$, which means that the collection of sets $\mathbb{S}_a=\{\phi\in\mathbb{S}(A):\phi(a)=1\textrm{ and }\phi(b)\geq\lambda\}$, for $a\in\mathcal{A}$, has the finite intersection property.  Again by the Banach-Alagolu Theorem, $\mathbb{S}=\bigcap_{a\in\mathcal{A}}\mathbb{S}_a$ is non-empty.  It is also convex and hence, by the Krein-Milman Theorem, contains extreme points.  We claim that any such extreme point of $\mathbb{S}$ is extreme in $\mathbb{S}(A)$, i.e. pure, which will complete the proof.  For if $\phi\in\mathbb{S}\backslash\mathbb{P}(A)$ then $\phi=\alpha\psi+(1-\alpha)\theta$ for some $\alpha\in(0,1)$ and $\psi,\theta\in\mathbb{S}(A)$.  As $1=\phi(a)=\alpha\psi(a)+(1-\alpha)\theta(a)$ and $\psi(a),\theta(a)\leq1$, we must in fact have $\psi(a)=\theta(a)=1$, for all $a\in\mathcal{A}$.  Thus $\psi(b)=\psi(a_1\ldots a_nba_n\ldots a_1)$, for all $n\in\omega$, which, as $||\psi||=1$ and $\inf\{||a_1\ldots a_nba_n\ldots a_1||:a_1,\ldots,a_n\in\mathcal{A}\}=\lambda$, implies that $\psi(b)\leq\lambda$.  Likewise, $\theta(b)\leq\lambda$ and hence, as $\lambda\leq\phi(b)=\alpha\psi(b)+(1-\alpha)\theta(b)$, we must in fact have $\psi(b)=\theta(b)=\lambda=\phi(b)$.  Thus $\psi$ and $\theta$ are both in $\mathbb{S}$ and hence $\phi$ is not extreme in $\mathbb{S}$.  The claim, and therefore the theorem, is thus proved. $\Box$\\

\begin{lem}\thlabel{phipsi}
Assume $A$ is a C$^*$-algebra and $\phi,\psi\in\mathbb{P}(A)$.  If $\phi^1_+\subseteq\psi^1_+$ then $\phi=\psi$.
\end{lem}

\paragraph{Proof:} Assume $\phi\neq\psi$.  As $\phi$ and $\psi$ are pure, the representations $\pi_\phi$ and $\pi_\psi$ are irreducible, by \cite{b} Theorem 1.6.6.  If $\pi_\phi$ and $\pi_\psi$ are inequivalent, we have $a\in A^1_+$ such that $\pi_\phi(a)v_\phi=v_\phi$ and $\pi_\psi(a)v_\psi=0$, by \cite{j} Corollary 7.  Then $\phi(a)=1$ and $\psi(a)=0$ so $a\in\phi^1_+\backslash\psi^1_+$ and we are done.  On the other hand, if $\pi_\phi$ and $\pi_\psi$ are equivalent, then we may assume they are in fact the same representation $\pi$ on a single Hilbert space containing both $v_\phi$ and $v_\psi$.  We can therefore write $v_\psi=v+\alpha v_\phi$ for some $\alpha\in\mathbb{F}$ and $v\perp v_\phi$.  By Kadison's Transitivity Theorem, we have $a\in A^1_+$ with $\pi(a)v=0$ and $\pi(a)v_\phi=v_\phi$.  Thus $\phi(a)=1$ and, as $\phi$ and $\psi$ are distinct and hence $v_\phi\neq\alpha v_\psi$, $|\psi(a)|\leq|\alpha|<1$, so again $a\in\phi^1_+\backslash\psi^1_+$. $\Box$\\

Pure states correspond to maximal proper left ideals (see \cite{i} Proposition 3.13.6(iv)).  Now we can show that they also correspond to maximal norm centred subsets.  We note that the formula (\ref{puremaxeq}) below actually follows from \cite{t} Proposition 2.2 where $\phi^1_+$, for $\phi\in\mathbb{P}(A)$, is even shown to excise $\phi$ in the sense that $\inf\{||aba-\phi(b)a^2||:a\in\phi^1_+\}=0$.

\begin{thm}\thlabel{puremax}
Assume $A$ is a C$^*$-algebra.  If $\phi\in\mathbb{P}(A)$ then $\phi^1_+$ is maximal norm centred and, for all $b\in A^1_+$,
\begin{equation}\label{puremaxeq}
\phi(b)=\inf\{||aba||:a\in\phi^1_+\}.
\end{equation}
If $\mathcal{A}\subseteq A^1_+$ is maximal norm centred then $\mathbb{S}(\mathcal{A})=\mathbb{P}(\mathcal{A})$ is a singleton.
\end{thm}

\paragraph{Proof:} Take $\phi\in\mathbb{P}(A)$ and extend $\phi^1_+$ to maximal norm centred $\mathcal{A}\subseteq A^1_+$.  For any $b\in A^1_+$, we have $\psi\in\mathbb{P}(A)$ with $\phi^1_+\subseteq\mathcal{A}\subseteq\psi^1_+$ and \[\psi(b)=\inf\{||a_n\ldots a_1\ldots a_nba_n\ldots a_1\ldots a_n||:n\in\omega\wedge a_1,\ldots,a_n\in\mathcal{A}\},\] by \thref{3.1(4)}.  But this means $\phi=\psi$, by \thref{phipsi}, and hence $\mathcal{A}=\phi^1_+$.

Now say $\mathcal{A}\subseteq A^1_+$ is maximal norm centred.  By \thref{3.1(4)}, $\mathbb{S}(\mathcal{A})=\{\phi\in\mathbb{S}(A):\mathcal{A}\subseteq\phi^1_+\}$ is non-empty.  Furthermore, for any $\phi\in\mathbb{S}(\mathcal{A})$ we actually have $\mathcal{A}=\phi^1_+$, as $\phi^1_+$ is norm centred and $\mathcal{A}$ is maximal norm centred.  By the Krein-Milman Theorem, $\mathbb{S}(\mathcal{A})$ is the closed convex hull of its extreme points, which must in fact be pure (see the argument at the end of the proof of \thref{3.1(4)}).  But for all $\phi,\psi\in\mathbb{P}(\mathcal{A})(=\mathbb{S}(\mathcal{A})\cap\mathbb{P}(A))$, we have $\phi^1_+=\mathcal{A}=\psi^1_+$ and hence $\phi=\psi$, by \thref{phipsi}. $\Box$\\

The following extension of the results in \cite{e} is required to prove \thref{KS}, specifically that the paving conjecture in \ref{KS4} implies the Kadison-Singer conjecture.  It can be thought of as analogous to the fact that, for any non-maximal filter $\mathcal{F}$ on $\omega$, there exists $X\subseteq\omega$ such that both $\{X\}\cup\mathcal{F}$ and $\{\omega\backslash X\}\cup\mathcal{F}$ generate a filter.

\begin{thm}\thlabel{Pperp}
Assume $A$ is a C$^*$-algebra, $\mathcal{A}\subseteq A^1_+$ is norm centred and $\mathcal{B}\subseteq A^1_+$ is a maximal norm centred extension of $\mathcal{A}$.  Then either $\mathcal{B}$ is the unique such extension or there exists $b\in\mathcal{B}$ such that $\mathcal{A}\cup\{1-b\}$ is norm centred.
\end{thm}

\paragraph{Proof:} Let $\psi\in\mathbb{P}(A)$ be such that $\psi^1_+=\mathcal{B}$. If $\mathcal{B}$ is not the unique maximal norm centred extension of $\mathcal{A}$, then we may take $a\in\mathcal{B}$ which is not in all of such extensions.  Thus $\inf\{||a_1\ldots a_n(1-a)a_n\ldots a_1||:a_1,\ldots,a_n\in\mathcal{A}\}=\epsilon>0$, for otherwise we would have $\phi(1-a)=0$ (see (\ref{phi1-a}) below) and hence $\phi(a)=1$, for all $\phi\in\mathbb{P}(\mathcal{A})$, and thus, by \thref{puremax}, $a$ would be in every maximal centred extension of $\mathcal{A}$, a contradiction.  By \thref{3.1(4)} we have $\phi\in\mathbb{P}(\mathcal{A})$ such that $\phi(1-a)=\epsilon$.  If $\psi$ and $\phi$ are inequivalent then we have $b\in A^1_+$ such that $\pi_\psi(b)v_\psi=v_\psi$ and $\pi_\phi(b)v_\phi=0$, by \cite{j} Corollary 7, and hence $b\in\mathcal{B}$ and $\phi(1-b)=1$, and we are done.  Otherwise, we may assume $H_\psi=H_\phi=H$ and $\pi_\psi=\pi_\phi=\pi$.  We claim that $v_\psi\perp v_\phi$.  To see this, take $a_1,\ldots,a_n\in\mathcal{A}$ such that $||\pi(c)||\leq||c||\leq\epsilon+\delta$, where $c=a_1\ldots a_n(1-a)a_n\ldots a_1$.  This implies that we have $\epsilon=\langle\pi(c)v_\phi,v_\phi\rangle\leq(\epsilon+\delta)\langle p^\perp v_\phi,v_\phi\rangle$, where $p$ denotes the projection onto the null space of $\pi(c)$.  But $\langle\pi(c)v_\psi,v_\psi\rangle=\psi(c)=\psi(1-a)=0$ so $pv_\psi=v_\psi$ and hence
\[|\langle v_\psi,v_\phi\rangle|^2=|\langle pv_\psi,v_\phi\rangle|^2=|\langle v_\psi,pv_\phi\rangle|^2\leq||pv_\phi||^2=\langle pv_\phi,v_\phi\rangle\leq1-\epsilon/(\epsilon+\delta)\rightarrow 0,\textrm{ as }\delta\rightarrow0.\]
Thus the claim is proved, and we have $b\in A^1_+$ such that $\pi(b)v_\psi=v_\psi$ and $\pi(b)v_\phi=0$, by the Kadison Transitivity Theorem, again giving $b\in\mathcal{B}$ and $\phi(1-b)=1$. $\Box$\\

\section{Norm Filters}\label{nfs}

While it is possible to define filters in $A^1_+$, for C$^*$-algebras $A$, in usual way (i.e. as directed upwards closed subsets), and we will indeed investigate these for projections in the Calkin algebra later on, these do not appear to be the most natural objects to study in the general C$^*$-algebra context.  In general, it appears to be the norm filters, as we define below, that are most relevant.

\begin{dfn}\thlabel{normfilter}
Assume $A$ is a C$^*$-algebra and $\mathcal{A}\subseteq A^1_+$.  We say $\mathcal{A}$ is a \emph{norm filter} if, for all $a\in A^1_+$, $\inf\{||a_1\ldots a_n(1-a)a_n\ldots a_1||:n\in\omega\textrm{ and }a_1,\ldots,a_n\in\mathcal{A}\}=0$ implies $a\in\mathcal{A}$.
\end{dfn}

We are keeping things symmetric here, but note that we could have equivalently defined $\mathcal{A}$ to be a norm filter if $\inf\{||(1-a)a_1\ldots a_n||:n\in\omega\textrm{ and }a_1,\ldots,a_n\in\mathcal{A}\}=0$ implies $a\in\mathcal{A}$, for all $a\in A^1_+$.  This is because $||b^*(1-a)b||\leq||(1-a)b||$ and $||(1-a)b||^2=||b^*(1-a)^2b||\leq||b^*(1-a)b||$, for all $a\in A^1_+$ and $b\in A$.

Also note that every proper norm filter is norm centred, for if we have $a_1,\ldots,a_n\in A^1_+$ such that $||a_1\ldots a_n||<1$ then $||(1-a)(a_1\ldots a_n)^m||\rightarrow0$, as $m\rightarrow\infty$, for all $a\in A^1_+$.

\begin{prp}\thlabel{Sinf}
If $A$ is a C$^*$-algebra and $\phi\in\mathbb{S}(A)$ then $\phi^1_+$ is a norm filter.
\end{prp}

\paragraph{Proof:}  Given $a\in A^1_+$ and $a_1,\ldots,a_n\in\phi^1_+$,
\begin{equation}\label{phi1-a}
\phi(1-a)=\phi(a_1\ldots a_n(1-a)a_n\ldots a_1)\leq||a_1\ldots a_n(1-a)a_n\ldots a_1||.
\end{equation}
Thus $\inf\{||a_1\ldots a_n(1-a)a_n\ldots a_1||:n\in\omega\textrm{ and }a_1,\ldots,a_n\in\phi^1_+\}=0$ implies $\phi(1-a)=0$ and hence $\phi(a)=1$, i.e. $a\in\phi^1_+$. $\Box$\\

\begin{cor}\thlabel{filtercentre}
Assume $A$ is a C$^*$-algebra and $\mathcal{A}\subseteq A^1_+$.  Then $\mathcal{A}$ is a (proper) norm filter if and only if it is a (non-empty) intersection of maximal norm centred subsets of $A^1_+$.
\end{cor}

\paragraph{Proof:} An intersection of norm filters is immediately verified to be a norm filter, so to prove the `if' part it suffices to show that maximal norm centred subsets of $A^1_+$ are norm filters.  By \thref{puremax}, this is equivalent to showing that $\phi^1_+$ is a norm filter, for all $\phi\in\mathbb{P}(A)$, which follows immediately from \thref{Sinf}.

On the other hand, say we have a norm filter $\mathcal{A}\subseteq A^1_+$.  For any $b\in A^1_+\backslash\mathcal{A}$, we have $\inf\{||a_1\ldots a_n(1-b)a_n\ldots a_1||:n\in\omega\textrm{ and }a_1,\ldots,a_n\in\phi^1_+\}>0$.  Thus, by \thref{3.1(4)}, there exists $\phi\in\mathbb{P}(\mathcal{A})$ with $\phi(1-b)>0$ and hence $\phi(b)<1$, i.e. $b\notin\phi^1_+$.  Therefore $\mathcal{A}=\bigcap_{\phi\in\mathbb{P}(\mathcal{A})}\phi^1_+$ which, by \thref{puremax}, is an intersection of maximal norm centred subset of $A^1_+$. $\Box$\\

In \cite{i} Theorem 1.5.2 it is shown that, for any C$^*$-algebra $A$, there are natural bijective order preserving correspondences between the collections of closed left ideals in $A$, hereditary C$^*$-subalgebras of $A$ and closed hereditary (real) cones in $A_+$.  In \cite{i} Theorem 3.10.7, it is further shown they correspond to weak$^*$ closed left invariant subspaces of $A^*$ and weak$^*$ closed faces of $\mathbb{S}(A)$.  The following corollary shows that they also correspond to norm filters in $A^1_+$.   

\begin{cor}
For any C$^*$-algebra $A$, the maps $\mathcal{A}\mapsto\{\lambda(1-a):\lambda\in\mathbb{R}^+\textrm{ and }a\in\mathcal{A}\}$ and $\mathcal{A}\mapsto\{1-a:a\in\mathcal{A}^1\}$ take norm filters in $A^1_+$ to closed hereditary cones in $A_+$ and vice-versa.
\end{cor}

\paragraph{Proof:} For any $\phi\in\mathbb{S}(A)$, the map $\mathcal{A}\mapsto\{1-a:a\in\mathcal{A}\}$ is immediately seen to take $\phi^1_+=\phi^{-1}[\{1\}]\cap A^1_+$ to $M^1_\phi=\phi^{-1}[\{0\}]\cap A^1_+$ and vice versa.  The image of any intersection of subsets of $\mathcal{A}^1_+$ is also immediately seen to be the intersection of the images under this map.  Thus the proof is complete by noting that any norm filter is an intersection of subsets of $A^1_+$ of the form $\phi^1_+$, by \thref{filtercentre}, and any closed hereditary cone is an intersection of subsets of $A_+$ of the form $M_\phi$, by \cite{i} Lemma 3.13.5 and Theorem 1.5.2. $\Box$\\

A subset of a lattice will be a filter if and only if it is upwards closed and closed under finite meets.  For arbitrary C$^*$-algebra $A$, $A^1_+$ may not be be a lattice, so we can not hope to prove exactly the same result.  However, we can use symmetric products in place of meets.  To this end, for use in the next proposition only, let us call $\mathcal{A}\subseteq A^1_+$, for C$^*$-algebra $A$, a \emph{product filter} if $\mathcal{A}$ is upwards closed (where $s\leq t\Leftrightarrow t-s\in A_+$) and $a_n\ldots a_2a_1a_2\ldots a_n\in\mathcal{A}$, whenever $n\in\omega$ and $a_1,\ldots,a_n\in\mathcal{A}$.  Also note in the following proof, and also later on in this article, we use the spectral family notation from \cite{r} so, for self-adjoint $s\in\mathcal{B}(H)$ and $t\in\mathbb{R}$, $E_s(t)$ is the spectral projection of $s$ corresponding to the interval $(-\infty,t]$.

\begin{prp}\thlabel{infprod}
Assume $A$ is a C$^*$-algebra and $\mathcal{A}\subseteq A^1_+$.  If $\mathcal{A}$ is a norm filter then it is a closed product filter.  If $A$ has real rank zero and $\mathcal{A}$ is a closed product filter then it is a norm filter.
\end{prp}

\paragraph{Proof:} By \thref{filtercentre}, proving the first part is equivalent to verifying that $\phi^1_+$ is a closed product filter for all $\phi\in\mathbb{P}(A)$, which is immediately seen to be true.

Conversely, assume $A$ has real rank zero and $\mathcal{A}$ is a closed product filter.  First we claim that if $a\in\mathcal{A}$, $p\in\mathcal{P}(A)$, $\epsilon>0$ and $p\geq E^\perp_a(1-\epsilon)$ then $p\in\mathcal{A}$.  To see this note that, for any $n\in\omega$ and $v\in\mathcal{R}(p)^\perp\subseteq\mathcal{R}(E_a(1-\epsilon))$ with $||v||=1$, we have $\langle a^nv,v\rangle\leq(1-\epsilon)^n$.  It follows that $a^n\leq p+(1-\epsilon)^np^\perp$.  As $\mathcal{A}$ is closed under symmetric products, $a^n\in\mathcal{A}$ and hence, as $\mathcal{A}$ is upwards closed, $p+(1-\epsilon)^np^\perp\in\mathcal{A}$.  But $p+(1-\epsilon)^np^\perp\rightarrow p$ as $n\rightarrow\infty$ and hence, as $\mathcal{A}$ is closed, $p\in\mathcal{A}$, which proves the claim.  Now say $a_1,\ldots,a_n\in\mathcal{A}$ and $a\in A^1_+$ satisfy $||a_1\ldots a_n(1-a)||<\epsilon$.  It follows that $||b(1-a)||<\epsilon$, where $b=a_n\ldots a_2a_1a_1\ldots a_n\in\mathcal{A}$.  As $A$ has real rank zero, there exists $p\in\mathcal{P}(A)$ such that $E_b^\perp(1-\epsilon/2)\leq p\leq E_b^\perp(1-\epsilon)$.  By the claim, $p\in\mathcal{A}$, and we also have $||p(1-b)||\leq\epsilon$.  Hence $||p(1-a)||\leq||p(1-b)(1-a)||+||pb(1-a)||\leq||p(1-b)||+||b(1-a)||\leq2\epsilon$.  Thus $p+p^\perp ap^\perp\in\mathcal{A}$ and $||p^\perp ap||\leq||p^\perp p||+||p^\perp(1-a)p||\leq||(1-a)p||\leq2\epsilon$, which yields $||a-(p+p^\perp ap^\perp)||\leq||a-pa-p^\perp ap^\perp||+||p-pa||=||p^\perp ap||+||p(1-a)||\leq4\epsilon$.  Thus if $\inf\{||a_1\ldots a_n(1-a)||:n\in\omega\textrm{ and }a_1,\ldots,a_n\in\mathcal{A}\}=0$ then, as $\mathcal{A}$ is closed, $a\in\mathcal{A}$ and hence $\mathcal{A}$ is a norm filter. $\Box$\\

\section{Projections}\label{op}

We already saw in \thref{infprod} how useful spectral projection approximations in C$^*$-algebras of real rank zero can be.  In fact, in C$^*$-algebras of real rank zero we can restrict our attention to just the projections $\mathcal{P}(A)$ in $A$, rather than the entirety of $A^1_+$, making the same definitions and proving the same theorems, with $\mathcal{P}(A)$ in place of $A^1_+$.  Indeed, this was the original approach in \cite{e}.

Specifically, the analog of \thref{phipsi} can be stated as follows.  If $A$ is a C$^*$-algebra of real rank zero, $\phi,\psi\in\mathbb{P}(A)$ and $\mathcal{P}(\phi)\subseteq\mathcal{P}(\psi)$ then $\phi=\psi$.  The proof is the same as the original, except that we replace $a$ with any projection $p$ such that $E^\perp_a(1-\epsilon)\leq p\leq E^\perp_a(\epsilon)$, for some $\epsilon\in(0,1/2)$.  This leads to obvious analogs of \thref{puremax} and \thref{Pperp}, although a little care has to be taken to prove the analog of the formula (\ref{puremaxeq}), namely $\phi(b)=\inf\{||pbp||:p\in\mathcal{P}(\phi)\}$.  To see this, note that if $a\in\phi^1_+$, $t\in(0,1)$ and $E^\perp_a(t+(1-t)/2)\leq p\leq E^\perp_a(t)$ then $p\in\mathcal{P}(\phi)$ (because $\pi_\phi(p)\geq E^\perp_{\pi_\phi(a)}(t+(1-t)/2)\geq E^\perp_{\pi_\phi(a)}(1-)$ and $E^\perp_{\pi_\phi(a)}(1-)v_\phi=v_\phi=\pi_\phi(a)v_\phi$) and $tp\leq a$ so $||pbp||=||b^{1/2}pb^{1/2}||\leq||b^{1/2}a^2b^{1/2}||/t^2=||aba||/t^2\rightarrow||aba||$ as $t\rightarrow1$.  In what follows we will refer to these theorems when their projection analogs are being used.


When restricting our attention to projections in the real rank zero case, we also have the following connection between norm filters that are filters and those that are countably closed (see \S\ref{order} for the definition of this and other standard order theoretic terminlogy). 

\begin{thm}\thlabel{sigfil}
Assume $A$ is a C$^*$-algebra and $\mathcal{P}$ is a norm filter in $\mathcal{P}(A)$.  If $\mathcal{P}$ is a filter and every $p\in\mathcal{P}(A)\backslash\{1\}$ is Murray-von Neumann below $p^\perp$, then $\mathcal{P}$ is countably closed.  Conversely, if $\mathcal{P}$ is countably closed and $A$ has real rank zero then $\mathcal{P}$ is a filter.
\end{thm}

\paragraph{Proof:} For the first part, take any strictly decreasing $(p_n)\subseteq\mathcal{P}\backslash\{1\}$ and let $u$ be a partial isometry such that $u^*u=p_0$ and $p_0u=0$.  Take any sequence $(\lambda_n)\subseteq(0,1/2)$ decreasing to $0$ and, for each $n\in\omega$, let $q_n$ be the projection onto $\mathcal{R}((1-\lambda_nu)(p_n-p_{n+1}))$ or, more precisely, let $q_n=E^\perp_{s_n}(0)$, where $s_n=(1-\lambda_nu)(p_n-p_{n+1})(1-\lambda_nu^*)$.  Note $q_n\in A$ by the functional calculus, as $\sigma(s_n)=\sigma((p_n-p_{n+1})(1-\lambda_n(u+u^*)+\lambda_n^2)(p_n-p_{n+1}))$ (because $\sigma(ab)\backslash\{0\}=\sigma(ba)\backslash\{0\}$ for arbitrary $a$ and $b$ in a C$^*$-algebra) and hence $\min(\sigma(s_n)\backslash\{0\})=1+\lambda_n^2-||\lambda_n(u+u^*)||>\lambda_n^2>0$.  Then the sequence $(\sum_{k\leq n}q_k+p_{n+1})\subseteq\mathcal{P}(A)$ is Cauchy and hence approaches some $p\in\mathcal{P}(A)$.  We then have $||p^\perp p_n||=\lambda_n/\sqrt{1+\lambda_n^2}\rightarrow0$, and hence $p\in\mathcal{P}$.  As $\mathcal{P}$ is a filter, we have $q\in\mathcal{P}$ with $q\leq p,p_0$ which, as $p_n=E^\perp_{p_0pp_0}(\frac{1-}{1+\lambda_n^2})$ for all $n\in\omega$, also satisfies $q\leq p_n$, for all $n\in\omega$, i.e. $\mathcal{P}$ is $\sigma$-closed.

On the other hand, say $A$ has real rank zero and take $p,q\in\mathcal{P}$.  Then, for any positive $(\lambda_n)$ with $\lambda_n\uparrow1$, we have $(p_n)\subseteq\mathcal{P}(A)$ with $E^\perp_{pqp}(\lambda_{n+1})\leq p_n\leq E^\perp_{pqp}(\lambda_n)$, for all $n\in\omega$.  It follows that $||p^\perp_n(pqp)^m||\leq\lambda_n^m\rightarrow0$, as $m\rightarrow\infty$, and hence $p_n\in\mathcal{P}$, for all $n\in\omega$.  As $\mathcal{P}$ is $\sigma$-closed, we have $r\in\mathcal{P}$ such that $r\leq p_n$, for all $n\in\omega$, and hence $r\leq p,q$. $\Box$\\

If we further restrict our attention to C$^*$-algebras $A$ of real rank zero such that $\mathcal{P}(A)$ is countably closed, then we see that the quantum objects we have defined can be described in purely order theoretic terms.

\begin{prp}\thlabel{nonquantum}
Assume $A$ is a C$^*$-algebra.  Then every centred $\mathcal{P}\subseteq\mathcal{P}(A)$ is a norm centred.  If $A$ has real rank zero and $\mathcal{P}(A)$ is countably closed then every norm centred $\mathcal{P}\subseteq\mathcal{P}(A)$ is centred.
\end{prp}

\paragraph{Proof:} This is essentially just \cite{a} Proposition 3.1, which we reprove here.  Given centred $\mathcal{P}\subseteq\mathcal{P}(A)$ and any $p_1,\ldots,p_n\in\mathcal{P}$, there exists $p\in\mathcal{P}(A)\backslash\{0\}$ with $p\leq p_1,\ldots,p_n$.  Then $p=pp_1\ldots p_n$ and hence $1=||p||=||pp_1\ldots p_n||\leq||p_1\ldots p_n||$, i.e. $\mathcal{P}$ is norm centred.

On the other hand, say $\mathcal{P}\subseteq\mathcal{P}(A)$ is norm centred, take $p_1,\ldots,p_n\in\mathcal{P}$ and set $a=p_n\ldots p_2p_1p_2\ldots p_n$, noting that $||a||=1$.  Also take $\lambda_m\uparrow1$ and $(q_m)\subseteq\mathcal{P}(A)$ such that, for all $m\in\omega$, $E^\perp_a(\lambda_{m+1})\leq q_m\leq E^\perp_a(\lambda_m)$, which is possible because $A$ has real rank zero.  As $||a||=1$, $q_m\neq0$ for any $m\in\omega$ so, as $\mathcal{P}(A)$ is countably closed, we have $q\in\mathcal{P}(A)\backslash\{0\}$ such that $q\leq q_m$, for all $m\in\omega$.  But, as $\lambda_m\uparrow1$, we must also have $q\leq p_n$ for all $n\in\omega$, i.e. $\mathcal{P}$ is centred. $\Box$\\

In these kinds of C$^*$-algebras, norm filters of projections can therefore be described in purely order theorectic terms as the intersections of maximal centred sets, by \thref{filtercentre}.  But beware that maximal centred sets are not necessarily filters, so norm filters are not necessarily filters either, even in these C$^*$-algebras.  A norm filter will, however, be the restriction of a filter in any lattice containing it, like its canonical Boolean completion for example.

\begin{prp}\thlabel{nflat}
If $A$ is a C$^*$-algebra of real rank zero, $\mathcal{P}(A)\backslash\{0\}$ is countably closed subset of a lattice $\mathbb{P}$, $\mathcal{P}$ is a norm filter in $\mathcal{P}(A)$ and $\mathcal{Q}$ is the filter generated by $\mathcal{P}$ in $\mathbb{P}$, then $\mathcal{P}=\mathcal{P}(A)\cap\mathcal{Q}$.
\end{prp}

\paragraph{Proof:} Take any $q\in\mathcal{P}(A)\cap\mathcal{Q}$.  As $\mathcal{Q}$ is the filter generated by $\mathcal{P}$, we have $p_1\wedge\ldots\wedge p_n\leq q$, for some $p_1,\ldots,p_n\in\mathcal{P}$.  This means that, for any $p\in\mathcal{P}$ with $p\leq p_1,\ldots,p_n$, we also have $p\leq q$.  Thus $q$ can be added to any centred subset containing $p_1,\ldots,p_n$ to form another centred subset.  This means that $q$ is in every maximal centred extension of $\mathcal{P}$ in $\mathcal{P}(A)$ and hence in $\mathcal{P}$. $\Box$\\

Unfortunately, we can not hope to get a converse of this, at least in general.  In fact, in any C$^*$-algebra $A$ satisfying the (first) hypotheses of \thref{sigfil}, consider the $(p_n)$ and $p$ in the proof.  Note that $p\notin\mathcal{P}$, where $\mathcal{P}$ is the p-filter given by the upwards closure of $(p_n)$, despite the fact that $p$ is in every norm filter containing $(p_n)$.  Thus $\mathcal{P}$ is not a norm filter even though it will be the restriction of a filter in any partial order containing it.

While norm filters will not necessarily be filters, the following proposition shows that in these C$^*$-algebras they will at least be closed under taking g.l.b.s.

\begin{prp}
Assume $A$ is a $C^*$-algebra of real rank zero, $\mathcal{P}(A)$ is countably closed and $p,q\in\mathcal{P}(A)$ have a g.l.b. $r\in\mathcal{P}(A)$.  Then any norm filter containing $p$ and $q$ will also contain $r$.
\end{prp}

\paragraph{Proof:} Any maximal centred subset containing $p$ and $q$ must also contain $r$, and hence the same is true of their intersections. $\Box$\\

We can also prove more order theoretic properties of certain subsets of projections in these C$^*$-algebras, like the following.

\begin{prp}
Assume $A$ is a $C^*$-algebra of real rank zero and $\mathcal{P}(A)$ is countably closed.  Then any centred $\mathcal{P}\subseteq\mathcal{P}(A)$ is in fact countably centred.
\end{prp}

\paragraph{Proof:} If some countable subset of $\mathcal{P}$ had no lower bound in $\mathcal{P}(A)\backslash\{0\}$ (i.e. if $0$ were its g.l.b.) then the same would be true of some finite subset, by \cite{a} Theorem 4.4, a contradiction. $\Box$\\

Note that a filter on an arbitrary preorder will be maximal centred (if and) only if it is maximal linked.  For projections in these C$^*$-algebras we do not even need the filter assumption, as shown in the following proposition.  In fact, the proposition and proof actually hold even if $\mathcal{P}(A)$ is not countably closed, so long as we replace `linked' with `norm linked' (defined analogously to norm centred).

\begin{prp}
Assume $A$ is a $C^*$-algebra of real rank zero and $\mathcal{P}(A)$ is countably closed.  Then any maximal centred $\mathcal{P}\subseteq\mathcal{P}(A)$ is maximal linked.
\end{prp}

\paragraph{Proof:} By \thref{nonquantum}, $\mathcal{P}$ is maximal norm centred, so by \thref{puremax}, we have $\phi\in\mathbb{P}(A)$ such that $\mathcal{P}=\mathcal{P}(\phi)$.  For any $q\in\mathcal{P}(A)\backslash\mathcal{P}$, we have $\phi(q)<1$ and hence there exist $p\in\mathcal{P}$ such that $||pq||^2=||pqp||<1$, by (\ref{puremaxeq}).  Thus $p$ and $q$ have no non-zero lower bound, by \cite{a} Proposition 3.1.  As $q$ was arbitrary, $\mathcal{P}$ is maximal linked. $\Box$\\

If $\mathcal{P}\subseteq\mathcal{P}(A)\backslash\{0\}$ has no non-zero lower bound then this means that, for every $q\in\mathcal{P}(A)\backslash\{0\}$ there exists $p\in\mathcal{P}(A)$ such that $q\nleq p$.  If $\mathcal{P}(A)$ were a Boolean algebra, this would mean $p^c$ and $q$ are compatible.  However, $\mathcal{P}(A)$ may well not be a Boolean algebra or even a lattice, and $q\nleq p$ is equivalent to the statement $||p^\perp q||>0$, while $p^\perp$ and $q$ being compatible is equivalent to the stronger statement that $||p^\perp q||=1$.  Nevertheless, for certain subsets $\mathcal{P}$ we can work a little harder and still obtain the stronger statement for some $p\in\mathcal{P}$.


\begin{prp}\thlabel{countdirect}
If $A$ has real rank zero and $\mathcal{P}(A)$ is countably closed then any countably directed $\mathcal{P}\subseteq\mathcal{P}(A)\backslash\{0\}$ will have no non-zero lower bound (if and) only if, for all $q\in\mathcal{P}(A)\backslash\{0\}$, there exists $p\in\mathcal{P}$ such that $p^\perp$ and $q$ are compatible.
\end{prp}

\paragraph{Proof:} Given $q\in\mathcal{P}(A)$, there always exists $p\in\mathcal{P}$ that is maximal for $||p^\perp q||$, as $\mathcal{P}$ is countably directed.  Take positive $(\lambda_n)$ with $\lambda_n\uparrow||p^\perp q||^2$.  As $A$ has real rank zero, we have $(p_n)\subseteq\mathcal{P}(A)$ such that $E^\perp_{qp^\perp q}(\lambda_{n+1})\leq p_n\leq E^\perp_{qp^\perp q}(\lambda_n)$, for all $n\in\omega$.  Then take $r\in\mathcal{P}(A)\backslash\{0\}$ with $r\leq p_n$, for all $n\in\omega$.  If $||p^\perp q||=1$ we are done, otherwise the projection $s(=E^\perp_{prp}(0))$ onto $\mathcal{R}(pr)$ is in $A$, by the functional calculus (note we are assuming here, as we may, that $A\subseteq\mathcal{B}(H)$ for some Hilbert space $H$).  As $\mathcal{P}$ has no lower bound then we can find $t\in\mathcal{P}$ such that $t\leq p$ and $s\nleq t$.  But then we have $v\in\mathcal{R}(r)$ such that $pv\notin\mathcal{R}(t)$ and hence $||p^\perp q||=||p^\perp v||<||t^\perp v||\leq||t^\perp q||$, contradicting the maximality of $||p^\perp q||$. $\Box$\\

\section{The Kadison-Singer Conjecture}\label{ksc}

The Kadison-Singer conjecture is a well-known long-standing conjecture stating that every pure state on an atomic MASA (maximal abelian subalgebra) of $\mathcal{B}(H)$, for infinite dimensional separable $H$, has a unique (pure) state exension.  We can use the theory we have developed so far to present some equivalent formulations.  Specifically, fix a basis $(e_n)$ for $H$ and, for $X\subseteq\omega$, define $P_X$ to be the projection onto $\overline{\mathrm{span}}\{x_n:n\in X\}$.  For $\mathcal{X}\subseteq\mathscr{P}(\omega)$ set $P_\mathcal{X}=\{P_X:X\in\mathcal{X}\}$ and let $\pi$ denote the canonical homomorphism from $\mathcal{B}(H)$ to the Calkin algebra $\mathcal{C}(H)=\mathcal{B}(H)/\mathcal{K}(H)$, where $\mathcal{K}(H)$ denotes the compact operators on $H$.  It is well known that $\mathcal{B}(H)$, and hence $\mathcal{C}(H)$ has real rank zero, and also that $\mathcal{P}(\mathcal{C}(H))$ is countably closed (see \cite{a}, for example), so we can indeed apply the theory developed so far.

\begin{thm}\thlabel{KS}
The following are equivalent.
\begin{enumerate}
\item\label{KS1} The Kadison-Singer conjecture.
\item\label{KS2} For every ultrafilter $\mathcal{U}\subseteq\mathscr{P}(\omega)$, $P_\mathcal{U}$ has a unique maximal norm centred extension.
\item\label{KS3} For every ultrafilter $\mathcal{U}\subseteq[\omega]^\omega$, $\pi[P_\mathcal{U}]$ has a unique maximal centred extension.
\item\label{KS5} For all $\epsilon>0$ and $P\in\mathcal{P}(\mathcal{B}(H))$ there exists $X_0,\ldots,X_{m-1}\subseteq\omega$ with $\bigcup_{k\in m}X_k=\omega$ such that $||PP_{X_k}||^2+||P^\perp P_{X_k}||^2\leq1+\epsilon$, for all $k<m$.

\item\label{KS4} There exists a $\delta>0$ such that, for all $P\in\mathcal{P}(\mathcal{B}(H))$ with $\langle Pe_n,e_n\rangle<\delta$, for all $n\in\omega$, there exists $X_0,\ldots,X_{m-1}\subseteq\omega$ with $\bigcup_{k\in m}X_k=\omega$ such that, for all $k\in m$, either $||PP_{X_k}||<1$ or $||P^\perp P_{X_k}||<1$.
\end{enumerate}
\end{thm}

\paragraph{Proof:}
\begin{itemize}
\item[\ref{KS1}$\Leftrightarrow$\ref{KS2}] By \thref{puremax}, any pure state $\phi$ on the atomic MASA $A=\{\sum\lambda_nP_{\{n\}}:(\lambda_n)\in l_\infty(\mathbb{F})\}$ is completely determined by $\mathcal{P}(\phi)$, which must be of the form $P_\mathcal{U}$ for an ultrafilter $\mathcal{U}$ on $\mathscr{P}(\omega)$, as $X\mapsto P_X$ is an isomorphism from $\mathscr{P}(\omega)$ onto $\mathcal{P}(A)$.  Likewise, any extension $\psi\in\mathbb{P}(\mathcal{B}(H))$ is determined by $\mathcal{P}(\psi)$, which must be a maximal norm centred subset of $\mathcal{P}(\mathcal{B}(H))$.

\item[\ref{KS1}$\Rightarrow$\ref{KS5}] Assume that $\epsilon>0$ and $P\in\mathcal{P}(\mathcal{B}(H))$ witness the failure of \ref{KS5}.  This means that \[\{X\subseteq\omega:||PP_X||^2+||P^\perp P_X||^2\leq 1+\epsilon\}\] generates a proper ideal of subsets of $\omega$, and hence there exists an ultrafilter $\mathcal{U}$ disjoint from it.  This, in turn, means that $\inf_{U\in\mathcal{U}}||P_UPP_U||+\inf_{U\in\mathcal{U}}||P_UP^\perp P_U||\geq 1+\epsilon$ and hence there exist necessarily distinct $\phi,\psi\in\mathbb{P}(\pi[P_\mathcal{U}])$ with $\phi(P)+\psi(P^\perp)\geq1+\epsilon$, by \thref{3.1(4)}.

\item[\ref{KS5}$\Rightarrow$\ref{KS4}] This follows immediately from the fact that $||PP_{X_k}||<1$ or $||P^\perp P_{X_k}||<1$ is equivalent to $||PP_{X_k}||^2+||P^\perp P_{X_k}||^2<2$.

\item[\ref{KS4}$\Rightarrow$\ref{KS2}] We assume that \ref{KS2} is false and prove that then \ref{KS4} must also be false.  So we have an ultrafilter $\mathcal{U}\subseteq\mathscr{P}(\omega)$ such that $P_\mathcal{U}$ does not have a unique maximal norm centred extension.  Define $\phi\in\mathbb{S}(\mathcal{B}(H))$ by $\phi(T)=\lim_{n\rightarrow\mathcal{U}}\langle Te_n,e_n\rangle$, and note that actually $\phi\in\mathbb{P}(\mathcal{B}(H))$ by \cite{l}.  Thus, by \thref{Pperp}, there exists $Q\in\mathcal{P}(\phi)$ such that $\{Q^\perp\}\cup P_\mathcal{U}$ is norm centred.  Given $\delta>0$, let $U\in\mathcal{U}$ be such that $\langle Qe_n,e_n\rangle>1-\delta/2$, for all $n\in U$, and let $P=E^\perp_{P_UQ^\perp P_U}(1/2)$.  As $\phi(P_U)=1$ and $\phi(Q)=1$, $\phi(P_UQ^\perp P_U)=0$ and hence $\phi(P)=0$.  But also, taking any $\psi\in\mathbb{S}(\{Q^\perp\}\cup P_\mathcal{U})$, we see that $\psi(P_UQ^\perp P_U)=1$ and hence $\psi(P)=1$.  Thus both $\{P\}\cup P_\mathcal{U}$ and $\{P^\perp\}\cup P_\mathcal{U}$ are norm centred.  In particular, for any $X_0,\ldots,X_{m-1}\subseteq\omega$ with $\bigcup_{k\in m}X_k=\omega$, there exists $k\in m$ such that $X_k\in\mathcal{U}$ and hence $||PP_{X_k}||=1=||P^\perp P_{X_k}||$.  But also, as $P\leq P_U$, we have $\langle Pe_n,e_n\rangle=0$, for all $n\in\omega\backslash U$, and $\langle Pe_n,e_n\rangle\leq\langle 2P_UQ^\perp P_Ue_n,e_n\rangle=\langle 2Q^\perp e_n,e_n\rangle\leq\delta$, for all $n\in U$.  Thus $P$ witnesses the failure of \ref{KS4}.

\item[\ref{KS2}$\Leftrightarrow$\ref{KS3}] If $\mathcal{U}\subseteq\mathscr{P}(\omega)$ is a principal ultrafilter, then $\mathcal{U}=\{P_U:n\in U\subseteq\omega\}$, for some $n\in\omega$, and hence $P_{\{e_n\}}\in P_\mathcal{U}$.  Then, for any $P\in\mathcal{B}(H)$ in some centred extension of $P_\mathcal{U}$, we have $1=||PP_{\{e_n\}}||=||Pe_n||$, which means $Pe_n=e_n$ and hence $P_{e_n}\leq P$, which means $P$ is in fact in every maximal centred extension of $P_\mathcal{U}$, i.e. this extension is unique.  Thus to verify \ref{KS2} it suffices to verify it only on the non-principal ultrafilters $\mathcal{U}$.  But if $\mathcal{V}$ is any norm centred extension of a non-principal ultrafilter $\mathcal{U}$, then $\pi[\mathcal{V}]$ will necessarily be (norm) centred.  Thus the maximal norm centred extensions of $P_\mathcal{U}$ and the maximal centred extensions of $\pi[P_\mathcal{U}]$ are in one-to-one correspondence. $\Box$\\
\end{itemize}

The statement \ref{KS5} above is essentially no different from paving conjectures already known to be equivalent to the Kadison-Singer conjecture.  It might be an appropriate version of the conjecture to apply, if it is indeed true.  On the other hand, \ref{KS4} above seems at first to be significantly weaker than \ref{KS5}, and would only be useful as a way of verifying the Kadison-Singer conjecture.  This appears to be new, although a Kadison-Singer equivalent statement quite close to \ref{KS4} is given in \cite{h} Conjecture 2.3, but for finite dimensional spaces.  However, there is a standard technique for turning such finite dimensional paving conjectures into infinite dimensional ones and vice versa, and the infinite dimensional equivalent of \cite{h} Conjecture 2.3 would be just like \ref{KS4} above, but with `$||PP_{X_k}||<1$ or $||P^\perp P_{X_k}||<1$' replaced by just `$||PP_{X_k}||<1$'.  In \cite{h}, the equivalence of Conjecture 2.3 to the Kadison-Singer conjecture is attributed to \cite{m} Theorem 1, and the proof there uses \cite{n} Proposition 7.7, which is somewhat similar to the proof of \ref{KS4}$\Rightarrow$\ref{KS2} given here.  The main difference is that, to obtain the $Q$ in the proof given here, we used the relatively elementary theory of states and norm centred sets, whereas the corresponding part of the proof of \cite{n} Proposition 7.7 uses the theory of supporting projections in the enveloping algebra of $\mathcal{B}(H)$, together with a noncommutive Urysohn lemma.

Once Kadison-Singer equivalent paving conjectures and their finite dimensional versions were discovered, most research on the Kadison-Singer conjecture seems to have become focused on these.  However, another approach, perhaps more in the spirit of the original formulation, is to try to prove the Kadison-Singer conjecture for certain kinds of states or, equivalently, certain kinds of (non-principal) ultrafilters on $\omega$.  Indeed, it might be that the conjecture fails in general and is only provable when certain extra assumptions are placed on the ultrafilters in question.

We will be interested in the following kinds of ultrafilters.

\begin{dfn}\thlabel{Pomega}
Order $[\omega]^\omega$ by $\subseteq^*$, where $A\subseteq^*B\Leftrightarrow|A\backslash B|<\infty$.  We say $\mathcal{U}\subseteq[\omega]^\omega$ is
\begin{enumerate}
\item a \emph{p-point} if $\mathcal{U}$ is a p-ultrafilter (i.e. a p-filter and an ultrafilter).
\item a \emph{q-point} if $\mathcal{U}$ is an ultrafilter and, for every interval partition $(I_n)$ of $\omega$, there exists $U\in\mathcal{U}$ such that $|U\cap I_n|\leq1$, for all $n\in\omega$.
\item \emph{rapid} if, for all $f\in\omega^\omega$, there exists $U\in\mathcal{U}$ such that $|U\cap f(n)|\leq n$, for all $n\in\omega$.
\end{enumerate}
\end{dfn}

The only previous result we know of in this direction is \thref{Reid}, which we now set up the necessary lemmas for proving.

\begin{lem}\thlabel{Reidlem}
For any $T\in\mathcal{B}(H)$, $f\in\omega^\omega$ and ultrafilter $\mathcal{U}\subseteq[\omega]^\omega$, there exists increasing $g\in\omega^\omega$ such that $g(n)\geq f(n)$, for all $n\in\omega$, and $\pi(P_UPP_U)=\pi(\sum_n P_{U\cap G_n}PP_{U\cap G_n})$, where $G_n=g(n+1)\backslash g(n)$, for all $n\in\omega$.
\end{lem}

\paragraph{Proof:} Set $g(0)=f(0)$ and, once $g(n)$ has been chosen, choose $g(n+1)\geq g(n),f(n+1)$ satisfying $||P_{g(n)}TP_{\omega\backslash g(n+1)}||\leq1/n^2$ and $||P_{\omega\backslash g(n+1)}TP_{g(n)}||\leq1/n^2$, which is possible because $P_{g(n)}T$ and $TP_{g(n)}$ have finite rank and are hence compact.  As $\sum1/n^2<\infty$, the operator $\sum P_{G_n}TP_{\omega\backslash g(n+2)}+P_{\omega\backslash g(n+2)}TP_{G_n}$ is compact.  As $\mathcal{U}$ is an ultrafilter, there is $U\in\mathcal{U}$ with $U\cap g(0)=0$ and, for all $n\in\omega$, $U\cap G_n=0$ or $U\cap G_{n+1}=0$.  It follows that, for all $n\in\omega$, $P_{U\cap G_n}TP_{U\cap G_{n+1}}=0=P_{U\cap G_{n+1}}TP_{U\cap G_n}$.  Thus
\begin{eqnarray*}
\pi(P_UTP_U) &=& \pi(\sum_nP_{G_n}TP_{\omega\backslash g(n+2)}+P_{\omega\backslash g(n+2)}TP_{G_n})\\
&& +\ \pi(\sum_nP_{U\cap G_n}TP_{U\cap G_{n+1}}+P_{U\cap G_{n+1}}TP_{U\cap G_n})\\
&& +\ \pi(\sum_n P_{U\cap G_n}TP_{U\cap G_n})\\
&=& \pi(\sum_n P_{U\cap G_n}TP_{U\cap G_n}).\Box\\
\end{eqnarray*}

\begin{lem}\thlabel{Reidlem2}
If $\mathcal{U}\subseteq[\omega]^\omega$ is an ultrafilter, $P\in\mathcal{P}(\mathcal{B}(H))$ and $\pi(P_UPP_U)=\pi(\sum_{n\in U}P_{\{n\}}PP_{\{n\}})$, for some $U\in\mathcal{U}$, then $\inf_{U\in\mathcal{U}}||P_UPP_U||=\lim_{n\rightarrow\mathcal{U}}\langle Pe_n,e_n\rangle=1-\inf_{U\in\mathcal{U}}||P_UP^\perp P_U||$.
\end{lem}

\paragraph{Proof:} The first equality follows from the fact that $||P_{\{n\}}PP_{\{n\}}||=\langle Pe_n,e_n\rangle$, for all $n\in\omega$.  But if $\pi(P_UPP_U)=\pi(\sum_{n\in U}P_{\{n\}}PP_{\{n\}})$ then \[\pi(P_UP^\perp P_U)=\pi(P_U-\sum_{n\in U}P_{\{n\}}PP_{\{n\}})=\pi(\sum_{n\in U}P_{\{n\}}-\sum_{n\in U}P_{\{n\}}PP_{\{n\}})=\pi(\sum_{n\in U}P_{\{n\}}P^\perp P_{\{n\}}),\] so $\inf_{U\in\mathcal{U}}||P_UP^\perp P_U||=\lim_{n\rightarrow\mathcal{U}}\langle P^\perp e_n,e_n\rangle=1-\lim_{n\rightarrow\mathcal{U}}\langle Pe_n,e_n\rangle$. $\Box$\\

\begin{thm}[Reid (1970)]\thlabel{Reid}
If $\mathcal{U}$ is a q-point then $\pi[P_\mathcal{U}]$ has a unique maximal centred extension.
\end{thm}

\paragraph{Proof:} Let $P\in\mathcal{P}(\mathcal{B}(H))$ be such that $\pi[P_\mathcal{U}]\cup\{\pi(P)\}$ is centred.  Take $U\in\mathcal{U}$ and $(G_n)$ be as in \thref{Reidlem} (with $f=0$).  As $\mathcal{U}$ is a q-point we may, by replacing $U$ with a subset if necessary, assume that $U\cap G_n$ contains at most one element, for all $n\in\omega$.  Thus $\pi(P_UPP_U)=\pi(\sum_{n\in U}P_{\{n\}}PP_{\{n\}})$ and hence $\inf_{U\in\mathcal{U}}||\pi(P_UP^\perp P_U)||\leq\inf_{U\in\mathcal{U}}||P_UP^\perp P_U||=0$, by \thref{Reidlem2}, i.e. $\pi(P)$ is in every maximal centred centred extension of $\pi[P_\mathcal{U}]$. $\Box$\\

We now show in \thref{rapp} that \thref{Reid} also holds for rapid p-points (giving an affirmative answer to the question raised in \cite{u} Problem 2 in the rapid case) instead of q-points and, furthermore, in this case the unique maximal centred extension is actually a filter.\footnote{Note that q-points are necessarily rapid but not vice-versa.  In fact, under MA(countable) there are rapid p-points that are not q-points - see \cite{g} Corollary 3 and Lemma 4.}

\begin{lem}\thlabel{UXmn}
If $\mathcal{U}\subseteq[\omega]^\omega$ is an ultrafilter, $m,n\in\omega$, $\mathcal{X}\subseteq[\omega]^{\leq m}$ and $|\{X\in\mathcal{X}:k\in X\}|\leq n$, for all $k\in\omega$, then there exists $U\in\mathcal{U}$ such that $|U\cap X|\leq1$, for all $X\in\mathcal{X}$.
\end{lem}

\paragraph{Proof:} For each $k\in\omega$, recursively choose $A_k$ to be a maximal subset of $\omega\backslash\bigcup_{j<k}A_j$ such that $|A_k\cap X|\leq 1$, for all $X\in\mathcal{X}$.  If $i\in\omega\backslash\bigcup_{j<k}A_j$, for some $i,k\in\omega$, then there exists $a_0,\ldots,a_{k-1}\in\omega$ and $X_0,\ldots,X_{k-1}\in\mathcal{X}$ such that $a_j\in A_j$ and $\{i,a_j\}\subseteq X_j$, for all $j<k$.  As there are at most $n(m-1)$ elements $a$ of $\omega$ such that $\{i,a\}\subseteq X$, for some $X\in\mathcal{X}$, we have must have $k\leq n(m-1)$, i.e. $\omega=\bigcup_{j\leq n(m-1)}A_j$ and hence $A_j\in\mathcal{U}$, for some $j\leq n(m-1)$. $\Box$\\

\begin{thm}\thlabel{rapp}
If $\mathcal{U}\subseteq[\omega]^\omega$ is a rapid p-point, $\pi[P_\mathcal{U}]$ is a maximal centred filter base.
\end{thm}

\paragraph{Proof:} Take any $P\in\mathcal{P}(\mathcal{B}(H))$ such that $\pi[P_\mathcal{U}]\cup\{\pi(P)\}$ is centred.  For each $n\in\omega$, let $\mathcal{X}_n=\{\{i,j\}:|\langle Pe_i,e_j\rangle|\geq n^{-3}\}$.  As $\sum_j|\langle Pe_i,e_j\rangle|^2\leq||Pe_i||^2\leq1$, we have $|\{X\in\mathcal{X}:i\in X\}|\leq n^6$, for all $i\in\omega$, and so we may apply \thref{UXmn} to get $U_n\in\mathcal{U}$ such that $|\langle Pe_i,e_j\rangle|<n^{-3}$, for all distinct $i,j\in U_n$.  As $\mathcal{U}$ is a p-point, we have $V\in\mathcal{U}$ and $f\in\omega^\omega$ such that $V\backslash f(n)\subseteq U_n$, for all $n\in\omega$.  Let $U$, $g$ and $(G_n)$ be as in \thref{Reidlem2} which, by replacing $U$ with $U\cap V$ if necessary, we may assume also satisfies $U\subseteq V$.  As $\mathcal{U}$ is rapid, we may further assume that $|U\cap g(n+1)|\leq n$, for all $n\in\omega$.  It follows that $||\sum_{i\neq j\in U\cap G_n}P_{\{i\}}PP_{\{j\}}||\leq n^2n^{-3}=n^{-1}$, for all $n\in\omega$, and hence $\pi(P_UPP_U)=\pi(\sum_{k\in U}P_{\{k\}}PP_{\{k\}})$.  By \thref{Reidlem2} $\inf_{U\in\mathcal{U}}||\pi(P_UP^\perp P_U)||\leq\inf_{U\in\mathcal{U}}||P_UP^\perp P_U||=0$, i.e. $\pi(P)$ is in every maximal centred centred extension of $\pi[P_\mathcal{U}]$. As $P$ was arbitrary, we have shown that the unique maximal centred extension of $\pi[P_\mathcal{U}]$ is $\mathcal{P}=\{p\in\mathcal{P}(\mathcal{C}(H)):\inf\{||p^\perp\pi(P_U)||:U\in\mathcal{U}\}=0\}$.  But as $\mathcal{U}$ is a p-point it follows that, for any $p\in\mathcal{P}$, the above infimum is actually attained by some $U\in\mathcal{U}$, i.e. we have $U\in\mathcal{U}$ such that $\pi(P_U)\leq p$. $\Box$\\

As rapidness was used in the proof of the above theorem simpy to show that $\pi[P_\mathcal{U}]$ has a unique maximal extension, if the Kadison-Singer conjecture is true then the above theorem holds for arbitrary p-point $\mathcal{U}$.  Conversely, if $\mathcal{U}$ is not a p-point, then no filter extension of $\pi[P_\mathcal{U}]$ can be a proper norm filter, as follows from the following proposition.

\begin{prp}\thlabel{nonp}
If $\mathcal{U}\subseteq[\omega]^\omega$ is a ultrafilter but not a p-point then there exists $p,q\in\mathcal{P}(\mathcal{C}(H))$ such that, for all $\phi\in\mathbb{S}(\pi[P_\mathcal{U}])$, $\phi(p)=1=\phi(q)$ even though $\phi(r)=0$, for all projections $r\leq p,q$.
\end{prp}

\paragraph{Proof:} Take decreasing $(X_n)\in\mathcal{U}$ with no pseudointersection in $\mathcal{U}$.  Let $p_n=\pi(P_n)$, for all $n\in\omega$, and take $p\in\mathcal{P}(\mathcal{C}(H))$ as in the proof of \thref{sigfil}, so $r\leq p_0,p$ if and only if $r\leq p_n$, for all $n\in\omega$.  But $r\leq p_n$, for all $n\in\omega$, (if and) only if $r\leq\pi(P_X)$, for some pseudointersection $X$ of the $(X_n)$, by \cite{f} Claim 2.5.10.  As $\mathcal{U}$ is an ultrafilter, it contains $\omega\backslash X$, for every such $X$, and hence $\phi(r)\leq\phi(\pi(P_X))=1-\phi(\pi(P_{\omega\backslash X}))=0$. $\Box$\\

So if the Kadison-Singer conjecture holds and $\mathcal{U}\subseteq\mathscr{P}(\omega)$ is an ultrafilter then the unique maximal centred extension of $\pi[P_\mathcal{U}]$ is a filter if and only if $\mathcal{U}$ is a p-point.  \thref{nonp} also raises the following question.

\begin{qst}
Does there exists a pure state $\phi$ on the Calkin algebra such that the follow holds?\footnote{the $p$ and $q$ here could equivalently be replaced with a decreasing sequence $(p_n)$, which might be an easier form of the problem to deal with.}
\begin{enumerate}
\item\label{AC1} For every $p,q\in\phi^1_+$ there exists $r\leq p,q$ with $\phi(r)>0$, and
\item\label{AC2} For some $p,q\in\phi^1_+$ we have $\phi(r)<1$, for all $r\leq p,q$.
\end{enumerate}
\end{qst}

If the Kadison-Singer conjecture holds, even just for p-points, then such a state would not be diagonalizable by any atomic MASA (because \ref{AC1} means it can not come from non-p-point ultrafilter on $\omega$, by \thref{nonp}, while \ref{AC2} means it can not come from a p-point, by the comment above), i.e. it would be a counterexample to Anderson's Conjecture.  Such counterexamples are known to (consistently) exist (see \cite{c} Theorem 6.46, for example), although it is not clear if these satisfy, or can be modified to satisfy, the above two conditions.

Next we show that, even if $\mathcal{U}$ is not a p-point, $\pi[P_\mathcal{U}]$ may still (consistently) be an ultrafilter base (note we are now talking about an ultrafilter base of projections, not subsets of $\omega$).

\begin{thm}\thlabel{qpufb}
If $\mathcal{U}\subseteq[\omega]^\omega$ is a q-point, $\pi[P_\mathcal{U}]$ is an ultrafilter base.
\end{thm}

\paragraph{Proof:} Take any $P\in\mathcal{P}(\mathcal{B}(H))$ and assume there exists a filter $\mathcal{P}\subseteq\mathcal{P}(\mathcal{C}(H))$ containing both $\pi[P_\mathcal{U}]$ and $\pi(P)$.  We show that $\pi(P)\geq\pi(P_U)$, for some $U\in\mathcal{U}$.  As in the proof of \thref{Reid}, we have $U\in\mathcal{U}$ such that $\pi(P_UPP_U)=\pi(D)$, where $D=\sum_{n\in U}P_{\{n\}}PP_{\{n\}}$.  As $\mathcal{P}$ is a filter, we have $p\in\mathcal{P}$ such that $p\leq\pi(P_U),\pi(P)$.  Take $(\lambda_n)\subseteq\mathbb{R}$ with $\lambda_n\uparrow1$ and note that, for each $n\in\omega$, $p\leq E^\perp_{\pi(D)}(\lambda_n)\leq\pi(E^\perp_D(\lambda_n))=\pi(P_{U_n})$, where $U_n=\{k\in U:\langle Pe_n,e_n\rangle>\lambda_n\}$.  But then there exists a pseudointersection $V$ of the $(U_n)$ such that $p\leq\pi(P_V)$, by \cite{f} Claim 2.5.10.  Thus $||\pi(P_{\omega\backslash V})p||=0$ which, as $p\in\mathcal{P}$ and hence $||P_UP||=1$, for all $U\in\mathcal{U}$, means that $\omega\backslash V\notin\mathcal{U}$ and hence $V\in\mathcal{U}$.  But, as $V$ is a pseudointersection of the $(U_n)$, $\pi(P_V)\leq\pi(P)$. $\Box$\\

Thus if $\mathcal{U}\subseteq[\omega]^\omega$ is a q-point, $\pi[P_\mathcal{U}]$ will have both a unique maximal centred extension and a unique ultrafilter extension, by \thref{Reid} and \thref{qpufb} respectively.  Note, however, that unless $\mathcal{U}$ is also a p-point, these extensions will be distinct, by \thref{nonp}, i.e. the unique maximal centred extension will properly contain the unique ultrafilter extension.  We also know that if $\mathcal{U}\subseteq[\omega]^\omega$ is a rapid p-point then $\pi[P_\mathcal{U}]$ has a unique maximal centred extension which is also the unique ultrafilter extension, by \thref{rapp}.  These results could be a mere coincidence, or they could perhaps point to some deeper connection.  For example, it might be that the Kadison-Singer conjecture fails in general but, for ultrafilter $\mathcal{U}\subseteq[\omega]^\omega$, $\pi[P_\mathcal{U}]$ has a unique maximal centred extension if and only if it has a unique ultrafilter extension.  At the very least, investigating the following Kadison-Singer analog for ultrafilters might well shed some light on the original Kadison-Singer conjecture, even though a positive or negative answer to this would not immediately appear to affirm or negate the original Kadison-Singer conjecture.

\begin{qst}
Does $\pi[P_\mathcal{U}]$ have a unique ultrafilter extension, for every ultrafilter $\mathcal{U}\subseteq[\omega]^\omega$?
\end{qst}

While we can not answer this question, we can at least show, in ZFC alone, that $\pi[P_\mathcal{U}]$ may not be an ultrafilter base, as it was in \thref{rapp} and \thref{qpufb}.

\begin{thm}
There are ultrafilters $\mathcal{U}\subseteq[\omega]^\omega$ for which $\pi[P_\mathcal{U}]$ is not an ultrafilter base.
\end{thm}
\paragraph{Proof:}  Take disjoint $(I_n)\subseteq[\omega]^{<\omega}$ such that $|I_n|\rightarrow\infty$ and let $\mathcal{U}$ be any ultrafilter extending the filter $\{X\subseteq\omega:|X\cap I_n|/|I_n|\rightarrow1\}$, so $\lim\sup|U\cap I_n|/|I_n|>0$, for all $U\in\mathcal{U}$.  Let $P$ be the projection onto $\overline{\mathrm{span}}\{\sum_{m\in I_n}e_m:n\in\omega\}^\perp$ and note that $\pi(P_U)\nleq\pi(P)$ even though $\mathcal{R}(P_U)\cap\mathcal{R}(P)$ is infinite dimensional, for all $U\in\mathcal{U}$.  So $\{\pi(Q):Q\in\mathcal{P}(\mathcal{B}(H))\wedge\exists U\in\mathcal{U}(\mathcal{R}(P_U)\cap\mathcal{R}(P)\subseteq\mathcal{R}(Q))\}$ is a filter properly containing the upwards closure of $\pi[P_\mathcal{U}]$. $\Box$\\

\section{Maximal Centred Filters}\label{mcf}

The existence of a rapid p-point, and even a selective ultrafilter (i.e. an ultrafilter that is simultaneously a p-point and a q-point, for which \thref{rapp} follows simply from Reid's result) is known to be consistent with ZFC \textendash\, eg. they can be added generically by forcing with $\mathscr{P}(\omega)/\mathrm{Fin}$, or constructed using CH.  It is also known that it is consistent they do not exist, in which case no maximal centred filter could come from extending $\pi[P_\mathcal{U}]$, for ultrafilter $\mathcal{U}$ on $\omega$, by \thref{nonp}.  However, there may exist states which are not diagonalized by any atomic MASA, so this does not necessarily mean there are no maximal centred filters.  To show this, we need to go back and analyze the model without p-points a little further, which we do in \thref{nupf}.

\begin{prp}\thlabel{phi0}
Assume $A$ is a (non-zero) unital C$^*$-algebra and $\mathcal{P}(A)$ is countably closed and has no atoms.  Then, for all $\phi\in\mathbb{S}(A)$, $\mathcal{P}(\phi)\neq\{1\}$.
\end{prp}

\paragraph{Proof:} Define $(P_n)\subseteq\mathcal{P}(A)$ by recursion as follows.  Let $P_0\in\mathcal{P}(A)$ be arbitrary and, once $P_n$ has been defined, let $Q_0,Q_1\in\mathcal{P}(A)\backslash\{0\}$ be such that $P_n=Q_0+Q_1$, which is possible because $\mathcal{P}(A)$ has no atoms.  Then $\phi(Q_0)+\phi(Q_1)=\phi(P)$ so $\phi(Q_k)\leq\phi(P)/2$ for some $k\in\{0,1\}$ and we may set $P_{n+1}=Q_k$.  Then we have $\phi(P_n)\leq2^n$ so, taking $P\in\mathcal{P}(A)\backslash\{0\}$ such that $P\leq P_n$, for all $n\in\omega$, which is possible because $\mathcal{P}(A)$ is countably closed, we see that $\phi(P)=0$.  Hence $\phi(P^\perp)=1$ and $P^\perp\neq1$. $\Box$\\

For $\phi\in\mathbb{S}(\mathcal{C}(H))$, define $\phi(X)=\phi(\pi(P_X))$.  Note that $\phi$ is then finitely additive and monotone w.r.t. $\subseteq^*$ on $\mathscr{P}(\omega)$.  Also note that, for $\epsilon>0$ and $(r_n)\subseteq\mathbb{R}_+$, with $r_0=0$ and $r_n\uparrow r\leq1$, we can recursively construct a subsequence, still with $r_0=0$, such that $\sum\sqrt{r_{n+1}-r_n}<\sqrt{r}+\epsilon$.

\begin{thm}\thlabel{C(H)pfil}
If $\phi\in\mathbb{P}(\mathcal{C}(H))$ and $\mathcal{P}(\phi)$ is a filter then $\mathcal{X}(\phi)=\{X\subseteq\omega:\phi(X)=1\}$ is a non-meagre p-filter.\footnote{It is open whether such filters exist in ZFC.  If there existed a model of ZFC without non-meagre p-filters then this theorem would immediately show that there are also no maximal centred filters of projections in the Calkin algebra in this model.  However, it is at least known that the non-existence of non-meagre p-filters has large cardinal consistency strength (see \cite{d} Corollary 4.4.15) so, even if such a model exists, our method here is still useful in showing that the non-existence of maximal centred filters has consistency strength equal to that of ZFC.}  In fact, for $(X_n)\subseteq\omega$, if $(X_n)$ is decreasing then it has a pseudointersection $X\subseteq\omega$ such that $\phi(X)=\inf\phi(X_n)$, while if $\phi(X_n)\not\rightarrow0$ then there exists increasing $(k_n)\subseteq\omega$ and $X\subseteq\omega$ such that $X\subseteq^*\bigcup_{j\in[k_n,k_{n+1})}X_j$, for all $n\in\omega$, and $\phi(X)>0$.
\end{thm}

\paragraph{Proof:} To see that $\mathcal{X}(\phi)$ is non-meagre, take any interval partition $(I_n)$ of $\omega$, let $B$ be the C$^*$-algebra generated by $(P_{I_n})$ and apply \thref{phi0} to $A=\pi[B]$.

Now say we have increasing $(X_n)\subseteq\omega$ with $X_0=0$.  By the comment above, given $\epsilon>0$ we may revert to a subsequence with $\sum\sqrt{\phi(X_{n+1}\backslash X_n)}<\sqrt{\sup\phi(X_n)}+\epsilon$.  As $\mathcal{P}(\phi)$ is a filter and therefore countably closed, by \thref{sigfil}, we have $p\in\mathcal{P}(\phi)$ such that $\phi(X_{n+1}\backslash X_n)=||\pi(P_{X_{n+1}\backslash X_n})p||^2$, for all $n\in\omega$, by \thref{puremax} (\ref{puremaxeq}).  Take any $P\in\mathcal{P}(\mathcal{B}(H))$ such that $\pi(P)=p$ and then define a sequence $(F_n)\subseteq[\omega]^{<\omega}$ such that $F_n\subseteq X_{n+1}\backslash X_n$, for each $n\in\omega$, and $\sum||P_{X_{n+1}\backslash(X_n\cup F_n)}P||<\sqrt{\sup\phi(X_n)}+\epsilon$.  Then, letting $Y=\bigcup X_n\backslash\bigcup F_n$, we see that $\sqrt{\phi(Y)}\leq||P_YP||<\sqrt{\sup\phi(X_n)}+\epsilon$ and $X_n\subseteq^*Y$, for all $n\in\omega$.  As $\epsilon>0$ was arbitrary, we can find decreasing $(Y_n)$ with $X_n\subseteq^*Y_n$, for all $n\in\omega$, and $\inf\phi(Y_n)=\sup\phi(X_n)$.  As $([\omega]^\omega,\subseteq^*)$ has no $(\omega,\omega)$-gaps, we can find $X\subseteq\omega$ such that $X_n\subseteq^*X\subseteq^*Y_n$, for all $n\in\omega$, and hence $\phi(Y)=\sup\phi(X_n)$.  Thus the equivalent statement for decreasing sequences stated in the theorem also holds, and the fact that $\mathcal{X}(\phi)$ is a p-filter follows from this and the fact that, whenever $X,Y\in\mathcal{X}(\phi)$ we have $\phi(\pi(P_{X\cap Y}))=\phi(\pi(P_XP_Y))=1$ and hence $X\cap Y\in\mathcal{X}(\phi)$.

For the last part, let $\epsilon=\lim\sup\phi(X_n)/2(>0)$.  Note if $\phi(X_{m_{i+1}}\backslash\bigcup_{j\in[k,m_i)}X_j)\geq\epsilon/3^{n+1}$ for $i<l$, where $k<m_0<\ldots<m_l$, then $l\epsilon\leq3^{n+1}$.  Thus we may choose $(k_n)$ so that, for $n\in\omega$, we have $\phi(X_j)>\epsilon$, for some $j\in[k_n,k_{n+1})$, and $\phi(X_m\backslash\bigcup_{j\in[k_n,k_{n+1})}X_j)<\epsilon/3^{n+1}$, for all $m\geq k_{n+1}$.  Then $\phi(\bigcap_{n\leq m}\bigcup_{j\in[k_n,k_{n+1})}X_j)>\epsilon-\sum_{n<m}\epsilon/3^{n+1}>\epsilon/2$, for all $m\in\omega$, and from the previous paragraph we thus have $X\subseteq^*\bigcup_{j\in[k_n,k_{n+1})}X_j$, for all $n\in\omega$, such that $\phi(X)\geq\epsilon/2>0$. $\Box$\\

\thref{C(H)pfil} can be interpreted as saying that such $\phi$ satisfy a weak version of normality.  Specifically, recall that, for Von Neumann algebra $A$, $\phi\in\mathbb{S}(A)$ is said to be \emph{normal} if, whenever $(a_\alpha)\subseteq A^1_+$ is a monotone decreasing net in $A^1_+$, $\phi(a_\alpha)\rightarrow\phi(\bigwedge_\alpha a_\alpha)$ (see \cite{i} 3.9.2).  If $A$ is an arbitrary C$^*$-algebra then $(a_\alpha)$ may not have a g.l.b., although we can extend the definition in this case simply by requiring that there exists some $a\in A$ below $(a_\alpha)$ with $\phi(a_\alpha)\rightarrow\phi(a)$.  Thus \thref{C(H)pfil} implies that $\phi\in\mathbb{P}(\mathcal{C}(H))$ is what might be termed `sequentially normal on projections' if (and only if, by \thref{sigfil}) $\mathcal{P}(\phi)$ is a filter.

\begin{thm}\thlabel{nupf}
If ZFC is consistent then it is also consistent with ZFC that there are no maximal centred filters in $\mathcal{P}(\mathcal{C}(H))$.
\end{thm}

\paragraph{Proof:} The statement holds in the model of ZFC without p-points constructed in \cite{d} \S4.4B.  This follows from the analog of \cite{d} Lemma 4.4.11 given below in \thref{biglem}.  First, we define the forcing notion in question.

\begin{dfn}
For $\mathcal{X}\subseteq\mathscr{P}(\omega)$, $\mathcal{P}(\mathcal{X})$ is the collection of functions $p$ such that $\mathrm{ran}(p)\subseteq\{0,1\}$ and $\mathrm{dom}(p)\in\mathcal{X}^c=\{\omega\backslash X:X\in\mathcal{X}\}$.
\end{dfn}

Note that a forcing notion $\mathbb{P}$ is said to be \emph{$\omega^\omega$-bounding} if, for every $p\in\mathbb{P}$ and every name $\dot{f}$ for a function in $\omega^\omega$, there exists $q\leq p$ and $g\in\omega^\omega$ such that $q\Vdash\forall n\in\omega(\dot{f}(n)\leq g(n))$.

\begin{lem}\thlabel{biglem}
If $\mathcal{X}$ is a non-meagre p-filter and $\dot{\mathbb{P}}$ is any $\mathcal{P}(\mathcal{X})^\omega$-name for an $\omega^\omega$-bounding forcing notion then $\mathds{1}\Vdash_{\mathcal{P}(\mathcal{X})^\omega*\dot{\mathbb{P}}}\forall\phi\in\mathbb{P}(\mathcal{C}(H))(\mathcal{P}(\phi)\textrm{ is a filter }\Rightarrow\mathcal{X}\nsubseteq\mathcal{X}(\phi))$.
\end{lem}

\paragraph{Proof:} We mimic the proof of \cite{d} Lemma 4.4.11.  Suppose we had some $\phi\in\mathbb{P}(\mathcal{C}(H))$ witnessing the failure of statement in $V[G][H]$ (where $G$ and $H$ are $\mathcal{P}(\mathcal{X})$ and $\dot{\mathbb{P}}$ generic sets respectively), i.e. such that $\phi\in\mathbb{P}(\mathcal{C}(H))$, $\mathcal{P}(\phi)$ is a filter and $\phi(X)=1$, for all $X\in\mathcal{X}$.  For the generic $G=(x_n)\in(2^\omega)^\omega$, we may assume, by switching the zeros and ones in $G$ (i.e. in each $x_n$) if necessary, that $\phi(X_n)\not\rightarrow0$, where $X_n=x_n^{-1}[\{1\}]$, for all $n\in\omega$.  Thanks to \thref{C(H)pfil} and the $\omega^\omega$-bounding property, we then have $((p_n),\dot{p})\in\mathcal{P}(\mathcal{X})^\omega*\dot{\mathbb{P}}$ and increasing $g,(k_n)\in V\cap\omega^\omega$, with $k_0=0$, such that $((p_n),\dot{p})\Vdash_{\mathcal{P}(\mathcal{X})^\omega*\dot{\mathbb{P}}}\phi(Y)>0$, where $Y=\bigcap_n g(n)\cup\bigcup_{j\in[k_n,k_{n+1})}X_j$.  We also have increasing $h\in V\cap\omega^\omega$ dominating $g$ everywhere such that $\bigcup_n\bigcup_{j\in[k_n,k_{n+1})}\mathrm{dom}(p_j)\backslash h(n)\in\mathcal{X}^c$.  For $n\in\omega$ and $j\in[k_n,k_{n+1})$ define $q_j(k)=p_j(k)$, for $k\in\mathrm{dom}(p_j)$, and $q_j(k)=0$, for $k\in[h(n),h(n+1))\backslash\mathrm{dom}(p_j)$ (and undefined elsewhere).  It then follows that $((q_n),\dot{p})\Vdash_{\mathcal{P}(\mathcal{X})^\omega*\dot{\mathbb{P}}}\forall n(\bigcup_{j\in[k_n,k_{n+1})}X_j\cap[h(n),h(n+1))\subseteq\bigcup_{j\in[k_n,k_{n+1})}\mathrm{dom}(p_j))$ and hence $((q_n),\dot{p})\Vdash_{\mathcal{P}(\mathcal{X})^\omega*\dot{\mathbb{P}}}Y\backslash h(0)\subseteq\bigcup_n\bigcup_{j\in[k_n,k_{n+1})}\mathrm{dom}(p_j)\backslash h(n)\in\mathcal{X}^c$.  But this gives us the contradiction $((q_n),\dot{p})\Vdash_{\mathcal{P}(\mathcal{X})^\omega*\dot{\mathbb{P}}}\phi(Y)=0$. $\Box$\\

This result, phrased in terms of states (and with a little extra theory), gives the following.

\begin{cor}\thlabel{nupfcor}
If ZFC is consistent then it is also consistent with ZFC that, for all $\phi\in\mathbb{P}(\mathcal{C}(H))$, there are $p,q\in\mathcal{P}(\mathcal{C}(H))$ such that $\phi(p)=\phi(q)=1$ even though $\max\{\phi(r):r\leq p,q\}<1$.
\end{cor}

\paragraph{Proof:} In the previous model, $\mathcal{P}(\phi)$ is not a filter for any $\phi\in\mathbb{P}(\mathcal{C}(H))$, and hence there are $p,q\in\mathcal{P}(\phi)$ such that $r\notin\mathcal{P}(\phi)$ for all $r\in\mathcal{P}(\mathcal{C}(H))$ such that $r\leq p,q$.  But, as (even non-linear) countable pregaps in $\mathcal{P}(\mathcal{C}(H))$ can always be (possibly not strictly) interpolated, by \cite{a} Theorem 4.6, $\{\phi(r):r\leq p,q\}$ has a maximum, which must therefore be less than $1$. $\Box$\\

\section{Towers}\label{towers}

Finally, we investigate towers (see \thref{po}\ref{tower}) of projections in C$^*$-algebras.  The only general result of note we have is the following corollary of \thref{countdirect}.

\begin{cor}
Assume $A$ has real rank zero and $\mathcal{P}(A)$ is countably closed.  If $\mathcal{P}\subseteq\mathcal{P}(A)$ is a tower in $\mathcal{P}(A)$ then, for all $q\in\mathcal{P}(A)$, there exists $p\in\mathcal{P}$ such that $||p^\perp q||=1$.
\end{cor}

\paragraph{Proof:} As $\mathcal{P}(A)$ is countably closed and $\mathcal{P}\subseteq\mathcal{P}(A)$ is a tower, it must also be countably closed and hence countably directed, so the result follows immediately from \thref{countdirect}. $\Box$\\

For the remainder of this section we investigate the question of whether $\pi[P_\mathcal{X}]$ will be a tower for certain towers $\mathcal{X}\subseteq[\omega]^\omega$ (again with respect to the $\subseteq^*$ order), which might be considered as a tower analog of the Kadison-Singer conjecture.  In \cite{o} Proposition 2.4 a tower $\mathcal{X}$ was constructed under CH for which $\pi[P_\mathcal{X}]$ is not a tower in $\mathcal{P}(\mathcal{C}(H))$.  We do the same under the weaker assumption of a certain cardinal invariant equality.  It follows from \thref{MAsigma} that some extra set theortic assumption like this is necessary, and that $\sup_n\mathfrak{m}(\sigma\textrm{-$n$-linked})$ could not be replaced with $\mathfrak{m}(\sigma\textrm{-centred})$.\footnote{the proof would still actually work with $\mathbf{non}(\mathcal{M})=\sup_n\mathfrak{m}(\sigma\textrm{-$n$-linked})$ replaced by the slightly weaker assumption that $\mathfrak{m}(\sigma\textrm{-$n$-linked})=\mathbf{cov}^*(\mathcal{ED}_\mathrm{fin})$ for some $n\in\omega$ (see \cite{s} Lemma 1.6.9).}

\begin{thm}\thlabel{nonM=m}
Let $(I_n)$ be a partition of $\omega$ into finite subsets such that $|I_n|\rightarrow\infty$ and let $P$ be the projection onto $\overline{\mathrm{span}}_{n\in\omega}(\sum_{k\in I_n}e_k)$.
If $\mathbf{non}(\mathcal{M})=\sup_n\mathfrak{m}(\sigma\textrm{-$n$-linked})$ then there is a tower $(X_\xi)_{\xi\in\mathfrak{t}}\subseteq[\omega]^\omega$ such that $\pi(P)\leq\pi(P_{X_\xi})$, for all $\xi\in\mathfrak{t}$.
\end{thm}

\paragraph{Proof:} First note that $\mathbf{non}(\mathcal{M})=\mathfrak{m}(\sigma\textrm{-$m$-linked})$, for some $m\in\omega$, as $\mathbf{non}(\mathcal{M})$ has uncountable cofinality.  Also, by the characterisation of $\mathbf{cov}(\mathcal{M})$ coming from \cite{d} Lemma 2.4.2, and the duality between $\mathbf{cov}(\mathcal{M})$ and $\mathbf{non}(\mathcal{M})$, there exists $(Y_\xi)_{\xi\in\mathbf{non}(\mathcal{M})}\subseteq[\omega]^\omega$ such that, for all $\xi\in\mathbf{non}(\mathcal{M})$ and $n\in\omega$, $|Y_\xi\cap I_n|=1$ and, for all $X\in[\omega]^\omega$, there exists $\xi\in\mathbf{non}(\mathcal{M})$ such that $Y_\xi\cap X$ is infinite.  Let $X_0=\omega$ and, for each $\xi\in\mathfrak{t}(=\mathbf{non}(\mathcal{M})$ because $\mathfrak{m}(\sigma\textrm{-$m$-linked})\leq\mathfrak{t}\leq\mathbf{non}(\mathcal{M}))$, let $X_{\xi+1}\subset^*X_\xi\backslash Y_\xi$ be such that $|X_{\xi+1}\cap I_n|/|I_n|\rightarrow1$, which is possible because $|X_\xi\cap I_n|/|I_n|\rightarrow1$.  If $\xi\in\mathfrak{t}$ is a limit ordinal then we construct $X_\xi$ as follows.  Let $\mathbb{P}$ be the partial order whose underlying set is the collection of $4$-tuples $(n,\Lambda,F,k)$ where $n,k\in\omega$, $\Lambda\in[\xi]^{<\omega}$, $F\in[\omega]^{<\omega}$ and $|I_j\cap\bigcap_{\zeta\in\Lambda}X_\zeta|/|I_j|\geq 1-1/(n+1)$, for all $j\geq k$, and where $(l,\Theta,G,j)\leq(n,\Lambda,F,k)$ if and only if $n\leq l$, $\Lambda\subseteq\Theta$, $k\leq j$, $F\subseteq G\subseteq\bigcap_{\zeta\in M}X_\zeta$, and $|I_i\cap G|/|I_i|\geq1-m/(n+1)$, for all $i\in j\backslash k$.  Given $n,k\in\omega$, $F\in[\omega]^{<\omega}$ and $\Lambda_0,\ldots,\Lambda_{m-1}\in[\xi]^{<\omega}$ such that $(n,\Lambda_j,F,k)\in\mathbb{P}$, for all $j\in m$, we may let $l\geq k$ be large enough that $|I_i\cap\bigcap_{\zeta\in\bigcup_{j\in m}\Lambda_j}X_\zeta|/|I_i|\geq 1-1/(n+1)$, for all $i\geq l$.  As $|I_i\cap\bigcap_{\zeta\in\Lambda_j}X_\zeta|/|I_i|\geq 1-1/(n+1)$, for all $i\geq k$, $|I_i\cap\bigcap_{\zeta\in\bigcup_{j\in m}\Lambda_j}X_\zeta|/|I_i|\geq 1-m/(n+1)$, for all $i\in l\backslash k$, and hence $p=(n,\bigcup_{j\in m}\Lambda_j,F\cup(\bigcup_{i\in l}I_i\cap\bigcap_{\zeta\in\bigcup_{j\in m}\Lambda_j}X_\zeta),l)\in\mathbb{P}$.  We also immediately see that $p\leq(n,\Lambda_j,F,k)$, for all $j\in m$.  Thus $\mathbb{P}_{n,F,k}=\{(n,\Lambda,F,k)\in\mathbb{P}:\Lambda\in[\xi]^{<\omega}\}$ is $m$-linked, for each $n,k\in\omega$ and $F\in[\omega]^{<\omega}$, and hence $\mathbb{P}$ is $\sigma$-$m$-linked.  It follows that there exists a filter $\mathcal{F}$ having non-empty intersection which all the dense sets $D_\zeta=\{(m,\Lambda,F,k)\in\mathbb{P}:\zeta\in\Lambda\}$, $E_n=\{(m,\Lambda,F,k)\in\mathbb{P}:m\geq n\}$ and $F_n=\{(m,\Lambda,F,k)\in\mathbb{P}:k\geq n\}$, for $\zeta\in\xi$ and $n\in\omega$.  Thus $X_\xi=\bigcup_{(m,\Lambda,F,k)\in\mathcal{F}}F$ satisfies $X_\xi\subseteq^*X_\zeta$, for all $\zeta<\xi$, and $|X_\xi\cap I_n|/|I_n|\rightarrow 1$.  This completes the recursion and we see that, for any $\xi\in\mathfrak{t}$, $||P_{I_n\backslash X_\xi}P||\rightarrow0$, because $|X_\xi\cap I_n|/|I_n|\rightarrow 1$, and hence $||\pi(P_{\omega\backslash X_\xi}P)||=0$, i.e. $\pi(P)\leq\pi(P_{X_\xi})$.  But also any $X\in[\omega]^\omega$ will have infinite intersection with $Y_\xi$, for some $\xi\in\mathfrak{t}$, and hence $X\nsubseteq^*X_{\xi+1}$ by our construction, i.e. $(X_\xi)$ is indeed a tower in $[\omega]^\omega$. $\Box$\\

It follows that, despite what might be assumed from Wofsey's construction using CH, the existence of towers $\mathcal{X}\subseteq[\omega]^\omega$ such that $\pi[P_\mathcal{X}]$ is not a tower does not actually have much to do with the size of the continuum.  Specifically, we can have such towers of length continuum, for arbitrary large continuum (under MA, for example), or we can have such towers of length $\aleph_1$, where $\aleph_1$ is strictly less than the continuum (in the Sacks model, for example).

Wofsey also showed in \cite{o} that there consistely exist towers $\mathcal{X}\subseteq[\omega]^\omega$ such that $\pi[P_\mathcal{X}]$ is a tower, specifically that this holds for towers that are added generically with finite conditions.  We take a different approach, using another cardinal equality (even weaker than the one in \thref{nonM=m}) to construct such a tower.  In this case, we do not know if this assumption is necessary, or whether there exists a better ZFC construction of such a tower (see \thref{ZFCntower}).

\begin{thm}\thlabel{b=t}
If $\mathfrak{b}=\mathfrak{t}$ then there is a tower $(X_\xi)_{\xi\in\mathfrak{t}}\subseteq[\omega]^\omega$ such that $(\pi(P_{X_\xi}))_{\xi\in\mathfrak{t}}$ is a tower.
\end{thm}

\paragraph{Proof:} By the proof of \cite{d} Lemma 4.4.12, $\mathfrak{b}=\mathfrak{t}$ implies that there exists a tower $(X_\xi)_{\xi\in\in\mathfrak{t}}\subseteq[\omega]^\omega$ such that, for any interval partition $(I_n)$ of $\omega$, there exists $\xi\in\mathfrak{t}$ such that $X_\xi$ is disjoint from $I_n$, for infinitely many $n\in\omega$.  Note that this means $\pi(P)\nleq\pi(P_{X_\xi})$ for any projection $P$ onto a block subspace with blocks in $(I_n)$ (i.e. such that there exists $(v_n)\subseteq H$ with $v_n\in\mathrm{span}_{k\in I_n}(e_k)$, for all $n\in\omega$, and $\mathcal{R}(P)=\overline{\mathrm{span}}(v_n)$).  Thus, to prove that $(X_\xi)_{\xi\in\in\mathfrak{t}}$ is the required tower, it suffices to show that the collection of $\pi(P)$, where $P$ is the projection onto a block subspace of $H$, is dense (in the order theoretic sense) in $\mathcal{P}(\mathcal{C}(H))$.  In order to prove this, simply take any projection $Q$ onto an infinite dimensional subspace and recursively pick increasing $n_k\in\omega$ and $v_k\in\mathcal{R}(Q)\cap\mathcal{R}(P_{\omega\backslash n_k})$ such that $||P_{\omega\backslash n_{k+1}}v_k||$ approaches $0$ fast enough that $\pi(P)=\pi(P')$, where $P$ and $P'$ are the projections onto $\overline{\mathrm{span}}(v_k)$ and $\overline{\mathrm{span}}(P_{n_{k+1}}v_k)$ respectively.  Then $\mathcal{R}(P')\subseteq\mathcal{R}(Q)$ and hence $\pi(P)=\pi(P')\leq\pi(Q)$, and also $P$ is the projection onto a block subspace in $H$ with blocks in $(I_k)$, where $I_k=[n_k,n_{k+1})$, for all $k\in\omega$.  As $Q$ was arbitrary, we are done. $\Box$\\

Thus we see that under MA or in the Sacks model, for example, we have both kinds of towers, namely towers $\mathcal{X}\subseteq[\omega]^\omega$ such that $\pi[P_{\mathcal{X}}]$ is a tower, as well as towers $\mathcal{X}\subseteq[\omega]^\omega$ such that $\pi[P_{\mathcal{X}}]$ is not tower.  We now show that consistently all towers $\mathcal{X}\subseteq[\omega]^\omega$ give rise to towers $\pi[P_{\mathcal{X}}]$.

\begin{lem}[\cite{q} Lemma 1]\thlabel{lsc}
Assume that $\mathbb{P}$ is a $\sigma$-centred forcing notion, $\kappa$ is an uncountable regular cardinal, $(X_\xi)_{\xi\in\kappa}\subseteq[\omega]^\omega$ is decreasing, $\phi$ is a finite lower semicontinuous submeasure on $\omega$ and $\Phi(X)=\lim_{n\rightarrow\infty}\phi(X\backslash n)$, for all $X\subseteq\omega$.  If $(X_\xi)_{\xi\in\kappa}$ has no pseudointersection $X$ with $\Phi(\omega\backslash X)<\Phi(\omega)$ then this remains true in any $\mathbb{P}$-generic extension.
\end{lem}

\paragraph{Proof:} Assume the lemma is false, so we have $p\in\mathbb{P}$, $k,n\in\omega$ and a $\mathbb{P}$-name $\dot{X}$ such that $p$ forces $\dot{X}$ to be a pseudointersection of the $(X_\xi)$ such that $\phi(\omega\backslash(\dot{X}\cup n))\leq\Phi(\omega)-1/k$.  Let $\mathbb{P}=\bigcup_{l\in\omega}\mathbb{P}_l$ and $Y_l=\{j\in\omega:\nexists q\in\mathbb{P}_l(q\leq p\wedge q\Vdash j\notin\dot{X})\}$, where $\mathbb{P}_l$ is centred, for all $l\in\omega$.

We claim that $\phi(\omega\backslash(Y_l\cup n))\leq\Phi(\omega)-1/k$, for all $l\in\omega$.  Otherwise, we would have $l,m\in\omega$ such that such that $\phi(m\backslash(Y_l\cup n))>\Phi(\omega)-1/k$.  Then, for each $j\in m\backslash(Y_l\cup n)$, we would have $p_j\in\mathbb{P}_l$ with $p_j\leq p$ and $p\Vdash m\notin\dot{X}$.  As $\mathbb{P}_l$ is centred, this means we would have a lower bound $q\leq p$ of the $p_j$ which therefore forces $m\backslash(Y_l\cup n)\subseteq\omega\backslash(\dot{X}\cup n)$ and hence $\phi(\omega\backslash(\dot{X}\cup n))>\Phi(\omega)-1/k$, a contradiction.

Thus $Y_l$ is not a pseudointersection of $(X_\xi)$, for any $l\in\omega$, an hence there exists $\xi\in\omega_1$ such that $Y_l\nsubseteq^*X_\xi$, for all $l\in\omega$.  We claim that $p\Vdash\dot{X}\nsubseteq^*X_\xi$.  To see this, take any $j,l\in\omega$ and $q\leq p$ with $q\in\mathbb{P}_l$.  There exists $i\geq j$ such that $i\in Y_l\backslash X_\xi$ which, by the definition of $Y_l$, means that $q\nVdash i\notin\dot{X}$.  Hence there exists $r\leq q$ such that $r\Vdash i\in\dot{X}$.  As $q$ and $j$ were arbitrary, the claim is proved, which contradicts the assumption that $p$ forces $\dot{X}$ to be a pseudointersection of the $(X_\xi)$. $\Box$\\

\begin{lem}\thlabel{preserve}
Assume $(\mathbb{P}_n,\dot{\mathbb{Q}}_n)$ is a finite support iteration of ccc forcings, $\kappa$ is an uncountable regular cardinal, $(X_\xi)_{\xi\in\kappa}\subseteq[\omega]^\omega$ is decreasing, $\phi$ is a finite lower semicontinuous submeasure on $\omega$ and $\Phi(X)=\lim_{n\rightarrow\infty}\phi(X\backslash n)$, for all $X\subseteq\omega$.  If $\mathbb{P}_n$ forces that $(X_\xi)_{\xi\in\kappa}$ has no pseudointersection $X$ with $\Phi(\omega\backslash X)<\Phi(\omega)$, for all $n\in\omega$, then this is also forced by $\mathbb{P}_\omega$.
\end{lem}

\paragraph{Proof:} Assume that the lemma is false, so we have a $\mathbb{P}_\omega$-name $\dot{X}$ and a $p\in\mathbb{P}_\omega$ forcing that $\dot{X}$ is a pseudointersection of $(X_\xi)_{\xi\in\kappa}$ and $\Phi(\omega\backslash\dot{X})<\Phi(\omega)$.  By reducing $p$ if necessary, we may also assume we have $j,k\in\omega$ such that $p$ forces that $\phi((\omega\backslash\dot{X})\backslash j)\leq\Phi(\omega)-1/k$.  Work in $V[G_n]$ for the moment, where $n>\max(\mathrm{supp}(p))$ and $G_n$ is any $\mathbb{P}_n$-generic containing $p$.  We have $(p^m_n)_{m\in\omega}\subseteq\mathbb{P}_\omega/\mathbb{P}_n$ and $Y_n\subseteq\omega$ such that $p^n_m\Vdash Y_n\cap m=\dot{X}\cap m$, for all $m\in\omega$.  In particular, we have $p^n_m\Vdash((\omega\backslash Y_n)\backslash j)\cap m=((\omega\backslash\dot{X})\backslash j)\cap m$ and $1\Vdash\phi(((\omega\backslash\dot{X})\backslash j)\cap m)\leq\phi((\omega\backslash\dot{X})\backslash j)\leq\Phi(\omega)-1/k$ and hence $\Phi(\omega\backslash Y_n)\leq\phi((\omega\backslash Y_n)\backslash j)<\Phi(\omega)$.  This means $Y_n$ is not a pseudointersection of $(X_\xi)_{\xi\in\kappa}$, and hence there exists some $\xi_n\in\kappa$ such that $Y_n\nsubseteq^*X_{\xi_n}$.

Now work in the ground model $V$, where we have names $(\dot{p}^n_m)$, $(\dot{Y}_n)$ and $(\dot{\xi}_n)$ for the $(p^n_m)$, $(Y_n)$ and $(\xi_n)$ respectively.  As each $\mathbb{P}_n$ is ccc and $\kappa$ has uncountable cofinality, we can find $\xi\in\kappa$ such that $p$ forces $\xi$ to be an upper bound for the $(\dot{\xi}_n)$, and hence forces $\dot{Y}_n\nsubseteq^*X_\xi$, for all $n$.

Take any $i\in\omega$ and $q\leq p$.  For any $n>\max(\mathrm{supp}(q))$, we can find $m\geq i$ and $r\in\mathbb{P}_n$ with $r\leq q$ such that $r\Vdash m-1\in\dot{Y}_n\backslash X_\xi$.  But then $r\hat{\,}\dot{p}^n_m\Vdash m-1\in\dot{X}\backslash X_\xi$.  Thus $p\Vdash\dot{X}\nsubseteq^*X_\xi$, a contradiction. $\Box$\\

These results were used in \cite{q} Theorem 2 to show that consistently there are no towers in the dual filter to any analytic p-ideal.  This argument, combined with \cite{f} Lemma 2.5.14 (originally a result of Juris Steprans), that the $X\subseteq\omega$ such that $\pi(P_X)$ is below some fixed projection is an analytic p-ideal, gives us the following theorem.  Rather than working with analytic p-ideals, however, we work directly with finite lower semicontinuous submeasures.

\begin{thm}\thlabel{MAsigma}
If ZFC is consistent then it is also consistent with ZFC and arbitrary large regular continuum that $\mathfrak{m}(\sigma\textrm{-centred})=\mathfrak{c}$ and $\pi[P_{\mathcal{X}}]$ is a tower, for every tower $\mathcal{X}\subseteq[\omega]^\omega$.
\end{thm}

\paragraph{Proof:} The model we construct is essentially the same as the standard model proving the consistency of MA+$\neg$CH (see \cite{k} Chapter VIII \S6, for example), except that we iterate with only $\sigma$-centred forcings, rather than all ccc forcings.  Specifically, start off with a model of ZFC where $\kappa$ is an uncountable regular cardinal satisfying $2^{<\kappa}=\kappa$ and do a $\kappa$-stage finite support iteratation of all $\sigma$-centred forcings of cardinality $<\kappa$ (which suffices to make MA hold for all $\sigma$-centred forcings, regardless of their cardinality, by \cite{k} Chapter II Lemma 3.1), each one iterated cofinally often and including those that appear in intermediate models of the iteration, which can be done by some book-keeping.  By genericity, it follows that $\mathfrak{m}(\sigma\textrm{-centred})=\mathfrak{c}$ in the final model.

On the other hand, say we had some tower $\mathcal{X}\subseteq[\omega]^\omega$ in the final extension such that $\pi[P_{\mathcal{X}}]$ is bounded below by $\pi(P)$, for some infinite rank $P\in\mathcal{P}(\mathcal{B}(H))$.  This means that $\Phi(\omega\backslash X)=0$, for all $X\in\mathcal{X}$, where $\phi$ is the finite lower semicontinuous submeasure given by $\phi(X)=||P_XP||$, for all $X\subseteq\omega$ (see \cite{f} Lemma 2.5.14).  As lower semicontinuous submeasures (and projections) are in natural correspondence with the reals, there exists some intermediate model $M_\xi$ containing $\phi$.  But $\mathcal{X}$ must have length $\kappa$, as otherwise it would appear at some earlier stage of the iteration, and then its Mathias forcing $\mathbb{M}(\mathcal{X})$ would also appear as a $\sigma$-centred forcing in the iteration, making $\mathcal{X}$ have pseudointersection in the next stage.  Thus there must exist some $X\in\mathcal{X}$ such that $Y\nsubseteq^*X$, for every $Y\in M_\xi$.  We can then pick $\zeta>\xi$ such that $X\in M_\zeta$.  Repeating this process $\omega_1$ times, we end up with an intermediate model $M_\alpha$ where $\mathrm{cf}(\alpha)=\omega_1$, and hence every real in $M_\alpha$ appears at some earlier stage in the iteration.  This means that $M_\alpha\cap\mathcal{X}$ is a tower in $M_\alpha$ (of cofinality $\aleph_1$), i.e. $M_\alpha\cap\mathcal{X}$ no pseudointersection $X\in M_\alpha$, and a fortiori none such that $\Phi(X)<\Phi(\omega)=1$.  However, we may repeat the process one more time to get $X\in\mathcal{X}$ which is a lower bound of $\mathcal{X}\cap M_\alpha$ and $\beta\in\kappa\backslash\alpha$ such that $X\in M_\beta$.  By \thref{lsc} and \thref{preserve}, this contradicts $\Phi(\omega\backslash X)=0$. $\Box$\\

The natural remaining question, to which we do not know the answer, is the following.

\begin{qst}\thlabel{ZFCntower}
Is it consistent with ZFC that $\pi[P_{\mathcal{X}}]$ is not a tower for any $\mathcal{X}\subseteq[\omega]^\omega$?
\end{qst}

\section{Appendix: Order Terminology}\label{order}

In this section we define the standard order theoretic terminology used throughout this article.

\begin{dfn}\thlabel{po}
Say we have a preorder $\mathbb{P}$ and $\mathbf{0}_\mathbb{P}=\{p\in\mathbb{P}:\forall q\in\mathbb{P}(p\leq q)\}$.  Then $\mathcal{F}\subseteq\mathbb{P}$ is
\begin{enumerate}
\item \emph{compatible} if $\mathcal{F}$ has a lower bound in $\mathbb{P}\backslash\mathbf{0}_\mathbb{P}$.
\item \emph{$n$-linked} if every $n$ element subset of $\mathcal{F}$ is compatible (\emph{linked} means 2-linked).
\item \emph{[countably] centred} if every finite [countable] subset of $\mathcal{F}$ is compatible.
\item \emph{$\sigma$-$n$-linked} [\emph{$\sigma$-centred}] if $\mathcal{F}$ is a countable union of $n$-linked [centred] subsets.
\item \emph{[countably] directed} if every finite [countable] subset of $\mathcal{F}$ has a lower bound in $\mathcal{F}$.
\item \emph{countably closed} if every decreasing sequence in $\mathcal{F}\backslash\mathbf{0}_\mathbb{P}$ has a lower bound in $\mathcal{F}\backslash\mathbf{0}_\mathbb{P}$.
\item a \emph{[p-]filter} if $\mathcal{F}$ is upwards closed and [countably] directed.
\item an \emph{ultrafilter} if it is a maximal proper filter.
\item a \emph{base} for $\mathcal{G}$ if every $q\in\mathcal{G}$ is above some $p\in\mathcal{F}$.
\item\label{tower} a \emph{tower} if $\mathcal{F}$ is a reverse well-ordered subset of $\mathbb{P}\backslash\mathbf{0}_\mathbb{P}$ with no lower bound in $\mathbb{P}\backslash\mathbf{0}_\mathbb{P}$.
\item \emph{dense} if, for every $p\in\mathbb{P}\backslash\mathbf{0}_\mathbb{P}$, there exists $q\in\mathcal{F}\backslash\mathbf{0}_\mathbb{P}$ below $p$.
\item \emph{predense} if the downwards closure of $\mathcal{F}$ is dense.
\end{enumerate}
\end{dfn}

\end{document}